\documentclass[12pt]{amsart}

\usepackage{amsmath}
\usepackage{amssymb}
\usepackage[all]{xy}

\numberwithin{equation}{section}

\newtheorem{thm}{Theorem}[section]
\newtheorem{lem}[thm]{Lemma}
\newtheorem{prop}[thm]{Proposition}

\newtheorem{rem}[thm]{Remark}

\newtheorem{exam-nota}[thm]{Example-Notation}

\newtheorem{dfn-nota}[thm]{Definition-Notation}

\newtheorem{dfn-lem}[thm]{Lemma-Definition}

\renewcommand{\qed}{\begin{flushright} {\bf Q.E.D.}\ \ \ \ \
                  \end{flushright} }
\newcommand{\beqa}{\begin{eqnarray*}}
\newcommand{\eeqa}{\end{eqnarray*}}

\newcommand{\id}{\mbox{${\rm id}$}}

\newcommand{\End}{\mbox{${\rm End}$}}
\newcommand{\Hom}{\mbox{${\rm Hom}$}}

\newcommand{\pr}{\mbox{${\rm pr}$}}

\newcommand{\fa}{\mbox{${\mathfrak a}$}}
\newcommand{\ft}{\mbox{${\mathfrak t}$}}

\newcommand{\fg}{\mbox{${\mathfrak g}$}}
\newcommand{\fq}{\mbox{${\mathfrak q}$}}
\newcommand{\fl}{\mbox{${\mathfrak l}$}}
\newcommand{\fs}{\mbox{${\mathfrak s}$}}

\newcommand{\fh}{\mbox{${\mathfrak h}$}}
\newcommand{\fn}{\mbox{${\mathfrak n}$}}

\newcommand{\fp}{\mbox{${\mathfrak p}$}}

\newcommand{\fb}{\mbox{${\mathfrak b}$}}

\newcommand{\fm}{\mbox{${\mathfrak m}$}}
\newcommand{\fu}{\mbox{${\mathfrak u}$}}

\newcommand{\cstar}{\mbox{${\C}^*$}}

\newcommand{\Gr}{{\rm Gr}}

\newcommand{\dw}{\dot{w}}
\newcommand{\dy}{\dot{y}}

\newcommand{\PR}{\mbox{${\mathbb P}$}}
\newcommand{\tildePR}{\tilde{\PR}}

\newcommand{\openp}{{\PR}_0}
\newcommand{\PVopen}{{\PR}_0(V)}
\newcommand{\PVstaropen}{{\PR}_0(V^*)}

\newcommand{\C}{\mbox{${\mathbb C}$}}
\newcommand{\Z}{\mbox{${\bf Z}$}}

\newcommand{\R}{\mbox{${\mathbb R}$}}

\begin{document}

\title[]{On the wonderful compactification}

\author[S. Evens]{Sam Evens}
\address{Department of Mathematics, University of Notre Dame, Notre Dame, 46556}
\email{sevens@nd.edu}

\author[B. Jones]{Benjamin F Jones}
\address{Department of Mathematics, University of Georgia, Athens, GA 30602-3704}
\email{bjones@math.uga.edu}

\begin{abstract}
These lecture notes explain the construction and basic properties
of the wonderful compactification of a complex semisimple group of adjoint
type. An appendix discusses the more general case of a semisimple symmetric
space.
\end{abstract}

\maketitle
\tableofcontents

\section{Introduction}
\label{sec_intro}

The purpose of these notes is to explain the construction of the
wonderful compactification of a complex semisimple group of adjoint
type. The wonderful compactification of a symmetric space was introduced
by DeConcini and Procesi \cite{DCP}, and has been extensively studied
in algebraic geometry. Of particular interest are the recent proofs
of the Manin conjecture for the compactification in \cite{STT} and
\cite{GMO}. Intuitively, the wonderful compactification gives information
about the group at infinity. In this connection, we mention the paper
by He and Thomsen \cite{HT} showing that closures of different
regular conjugacy classes coincide at infinity. The wonderful compactification
has been used by Ginzburg 
 in the study of character sheaves \cite{Gi},
and by Lusztig in his study of generalized character sheaves \cite{Lu}.
It is closely related to the Satake compactification and to various
analytic compactifications \cite{BorelJi}.
It also plays an important role in Poisson geometry, as it is crucial for
undersanding the geometry of a moduli space of Poisson
homogeneous spaces \cite{EL1, EL2}. This list is not meant to be exhaustive,
but to give some idea of how the wonderful compactification is
related to the rest of mathematics. The reader may consult \cite{springericm}
and \cite{timashev} for further references.

The appendix contains a construction of the wonderful compactification
of a complex symmetric space due to DeConcini and Springer \cite{DCS1}.
This construction does not specialize to the construction in 
Section \ref{sec_construction} in the case when the symmetric space
is a group, but it is equivalent, and the two constructions are
conceptually similar.

These notes are based on lectures given by the first author at
Hong Kong University of Science and Technology and by both authors
at Notre Dame. They are largely based on work of DeConcini, Procesi,
and Springer which is explained in the papers \cite{DCP},
\cite{DCS1}, and the book \cite{BK}. The notes add little in content,
and their purpose is to make the simple construction of the wonderful
compactification more accessible to students and to mathematicians
without an extensive background in algebraic groups and algebraic
geometry.

We would like to thank Hong Kong University of Science and Technology
for hospitality during preparation of a first draft of these notes,
and would like to thank Jiang-Hua, Dragan Milicic, Allen Moy, Xuhua He,
Francois Ledrappier, and Dennis Snow for useful comments and discussions. 

\subsection{Notation}\label{intro-notation}
In these notes, an algebraic variety is a complex quasi-projective
variety, not necessarily irreducible. A subvariety is a locally closed
subset of a variety.

\section{Construction and basic properties of the wonderful compactification}
\label{sec_construction}

We explain how to prove that the wondeful compactification
of a semisimple complex group $G$ of adjoint type is
smooth, and describe its $G\times G$-orbit structure.

\subsection{Definition of the compactification}
\label{sec_definition}

Let $G$ be a complex connected semisimple group with trivial center.
Let $\tilde{G}$ be the simply connected cover of $G$, and
 choose a maximal torus $\tilde{T}$ contained in a Borel subgroup $\tilde{B}$
of $\tilde{G}$. We denote their images in $G$ by $T\subset B$. We denote Lie
algebras of algebraic groups by the corresponding gothic letter, so the 
Lie algebras of $T \subset B \subset G$, and of $\tilde{T} \subset \tilde{B}
\subset \tilde{G}$, are $\ft \subset\fb \subset \fg$.
Let $X^*(\tilde{T})$ be the group of characters of 
$\tilde{T}$. We write this group additively, so if
$\lambda, \mu \in X^*(\tilde{T})$ we have by definition 
$(\lambda+\mu)(t) = \lambda(t) \mu(t)$ for all $t \in \tilde{T}$. 
Let $\Phi \subset X^*(\tilde{T})$ be the roots of $\tilde{T}$ in $\tilde{G}$,
and take the positive roots $\Phi^+$ to be the roots of $\tilde{T}$ in
$\tilde{B}$. Let $\{ \alpha_1, \dots, \alpha_l \}$ be the corresponding
set of simple roots, where $l=\dim(\tilde{T})$.
With these choices, for $\lambda, \mu \in X^*(\tilde{T})$, we say $\lambda 
\ge \mu$  if $\lambda - \mu = \sum_{\alpha \in \Phi^+} n_\alpha \alpha$ where
 the $n_\alpha$ are nonnegative integers. There is an embedding of the
character group $X^*(T) \hookrightarrow X^*(\tilde{T})$ as the characters
that are trivial on the center of $G$. If $\lambda$ is in the image of this
embedding, we may compute $\lambda(t)$ for $t\in T$.

Let $\tilde{B}^-$ be the opposite Borel of $\tilde{G}$ containing $\tilde{T}$
and let $B^-$ be its image in $G$.
Let $U$ and  $U^-$ be the unipotent radicals of $\tilde{B}$ and $\tilde{B}^-$.
The group homomorphism $\tilde{G} \to G$ restricts to an isomorphism on any unipotent
subgroup, and we identify $U$ and $U^-$ with their images in $G$.
If $W$ is a representation of $\tilde{G}$, set $W_\mu$ for the
$\mu$-weight space for the $\tilde{T}$-action.  If $v\in W_\mu$,
then $U^- \cdot v \subset v + \sum_{\phi < \mu} W_\phi$.

Fix an irreducible representation $V=V(\lambda)$ of $G$ with regular 
highest weight $\lambda$ and choose a basis 
$v_0, \dots, v_n$ of $T$-weight vectors of  $V$ with the following properties:

(1) $v_0$ has weight $\lambda$;

(2) For $i \in \{ 1, \dots, l \},$ $v_i$ has weight 
$\lambda - \alpha_i$;

(3) Let $\lambda_i$ be the weight of $v_i$. Then if
$\lambda_i < \lambda_j$, then $i > j$.

\begin{rem}
\label{rem_weightmultiplicity}
{\em
For $i=0, \dots, l$, $\dim(V_{\lambda_i}) = 1$.
}
\end{rem}

\begin{rem}
\label{rem_lowertriangular}
{\em
$U^- \cdot v_i - v_i$ is in the span of the $v_j$ for
$j>i$.
}
\end{rem}

The induced action of $\tilde{G}$ on the projective space $\PR (V)$ factors
to give an action of $G$. 
Define $\PVopen = \{ [ \sum_{i=0}^{n} b_i v_i ] : b_0 \not= 0 \} \subset
 \PR (V)$.
The affine open set $\PVopen \cong  \{ v_0 + \sum_{i=1}^{n} b_i v_i \} \cong \C^n$. 
Further, $\PVopen$ is $U^-$-stable by Remark \ref{rem_lowertriangular}
applied to $v_0$.

Also, the morphism $U^- \to U^-\cdot [v_0]$ is an isomorphism of
algebraic varieties. Indeed,  the stabilizer of $[v_0]$ in
$G$ is $B$, so the stabilizer in $U^-$ of $[v_0]$ is trivial, and 
the result follows.

\begin{lem}
\label{lem_uminusclosed}
$U^- \mapsto U^-\cdot [v_0]$ is an isomorphism between
$U^-$ and the closed subvariety $U^- \cdot [v_0]$ of
$\PVopen$.
\end{lem}

We have proved everything but the last statement, which
follows by the following standard fact.

\begin{lem}
\label{lem_unipotentorbitclosed}
Let $A$ be a unipotent algebraic group and let $X$ be an affine
$A$-variety. Then every $A$-orbit in $X$ is closed.
\end{lem}

$\tilde{G}$ acts on $V^*$ by the dual action, $g\cdot f = f\circ g^{-1}$.
Choose a dual basis $v_0^*, \dots, v_n^*$ of $V^*$. Then
each $v_i^*$ is a weight vector for $\tilde{T}$ of weight $-\lambda_i$,
and $v_0^*$ is a highest weight vector with respect to
the negative of the above choice of positive roots. 
Note that $\PVopen = \{ [v] : v_0^*(v)\not= 0 \}$.

Define $\PVstaropen = \{ [\sum_{i= 0}^n b_i v_i^* ] : b_0 \not= 0 \}$.
Then $U$ acts on $\PVstaropen$, and the morphism $U\mapsto U\cdot [v_0^*]$
is an isomorphism of varieties.  We state the analogue
of Lemma \ref{lem_uminusclosed} for this action, which follows 
by reversing the choice of positive roots.
\begin{lem}
\label{lem_uclosed}
$U \mapsto U\cdot [v_0^*]$ is an isomorphism between
$U$ and the closed subvariety $U \cdot [v_0^*]$ of
$\PVstaropen$.
\end{lem}

$\tilde{G} \times \tilde{G}$ acts on $\End(V)$ by the formula
$(g_1, g_2)\cdot A = g_1 A g_2^{-1}$ and on $V\otimes V^*$ by linear extension of the
formula
$(g_1, g_2)\cdot v\otimes f = g_1\cdot v \otimes g_2\cdot f$.
Then the canonical identification $V\otimes V^* \cong \End(V)$
is $\tilde{G}\times \tilde{G}$-equivariant, and we treat this
identification as an equality. 
 Then $\{ v_i \otimes v_j^* : i, j = 0, \dots n\}$ is a basis of $\End(V)$.
 As before, the $\tilde{G} \times \tilde{G}$-action
on $\PR (\End(V))$ descends to give a $G\times G$-action.

 We consider the open set 
 defined by
$$\openp = \{ [ \sum a_{ij}v_i \otimes v_j^*] : a_{00} \not= 0 \}
      =\{ [A] \in \PR (\End(V)) : v_0^* (A\cdot v_0) \not= 0 \}.$$
It follows from the above that $U^-T \times U$ preserves
$\openp$. Indeed, we use that fact that $U^-T$ preserves 
$\PVopen$ and $U$ preserves  $\PVstaropen$.

Define an embedding $\psi:G\to \PR (\End(V))$ by
$\psi(g)=[g]$, and note that $\psi$ is $G\times G$-equivariant, 
where $G\times G$ acts on $G$ by $(g_1, g_2)\cdot x = g_1 x g_2^{-1}$.
Let $X=\overline{\psi(G)}$. 

\subsection{Geometry of the open affine piece}
\label{sec_structureopenaffine}

Let $X_0 = X\cap \openp$.
Let the Weyl group $W=N_G(T)/T$, and for each $w\in W$, choose a representative
$\dw \in N_G(T)$.

\begin{lem}
\label{lem_X0G} 
 $X_0 \cap \psi(G) = \psi(U^-TU)$.
\end{lem}

\noindent
{\bf Proof.} Clearly, $\psi(e)\in X_0$, and
$\psi(U^-TU)=(U^-T \times U)\cdot \psi(e)$,
so 
$\psi(U^-TU) \subset X_0$ since $\openp$ is
$U^-T \times U$ stable. This gives one inclusion.

For the other inclusion, recall that by the Bruhat
decomposition, 
$G = \cup_{w\in W} U^-T w U$. Then 
 $\dw\cdot [v_0]$ is a weight
vector of weight $w\lambda$. Since $\lambda$ is regular, if $w\not= e$, then $\dw \cdot [v_0] \not\in X_0$. Thus, $\psi(\dw)\not\in
X_0$. Now use $U^-T \times U$-stability of $X_0$
again.
\qed

Since $\psi(U^-TU)$ is dense in $\psi(G)$, it is dense in $X$, and it follows that
$X_0$ is the closure  $\overline{\psi(U^-TU)}$ in $\openp$.

Define $Z=\overline{\psi(T)}$, the closure of $\psi(T)$
in $\openp$.

\begin{lem}
\label{lem_Zaffine}
$Z \cong {\C}^l$, where $l=\dim(T)$.
\end{lem}

\noindent
{\bf Proof.} Under the identification of 
$\End(V)$ with $V\otimes V^*$, the identity corresponds
to $\sum_{i=0}^n v_i \otimes v_i^*$. Thus, for $\tilde{t}\in \tilde{T}$ projecting
to $t\in T$,
$$\psi(t) = t\cdot \psi(e) = [\sum \tilde{t}\cdot v_k \otimes v_k^* ]
=[\sum \lambda_k(\tilde{t}) v_k \otimes v_k^*]$$
This is
$$
[\lambda(\tilde{t})v_0 \otimes v_0^* + \sum_{k > 0} \lambda_k(\tilde{t})v_k \otimes
v_k^*] = 
[v_0 \otimes v_0^* + \sum_{k > 0} \frac{\lambda_k(\tilde{t})}{\lambda(\tilde{t})}
v_k \otimes v_k^*].$$

\noindent
But 
$$ 
\frac{\lambda_k(\tilde{t})}{\lambda(\tilde{t})}=\frac{1}{\prod_{i=1}^{l}\alpha_i(t)^{n_{ik}}}
$$
where $\lambda_k = \lambda - \sum_{i=1}^{l} n_{ik}\alpha_i$.

Thus, 
$$\psi(t) = 
[v_0 \otimes v_0^* + \sum_{i=1}^l \frac{1}{\alpha_i(t)}v_i\otimes v_i^*
+ \sum_{k>l} \frac{1}{\prod_{i=1}^{l}\alpha_i(t)^{n_{ik}}} v_k \otimes v_k^* ]
.$$

Define $F:{\C}^l \to \openp$ by
$$
F(z_1, \dots, z_l) = 
[v_0 \otimes v_0^* + \sum_{i=1}^l z_i v_i\otimes v_i^*
+ \sum_{k>l} (\prod_{i=1}^{l} z_i^{n_{ik}}) v_k \otimes v_k^* ]
.$$
Then $F$ is a closed embedding, so to prove the Lemma, it suffices to identify the
image with $\overline{\psi(T)}$. This is routine since by the above calculation,
$\psi(T)$ is in the irreducible $l$-dimensional variety given by the image of $F$.
\qed

Note that $U^- \times U$ stabilizes $X_0 = X\cap \openp$, and
define $\chi:U^- \times U \times Z \to X_0$ by
the formula $\chi(u,v,z)=(u,v)\cdot z = uzv^{-1}$.

\begin{thm}
\label{thm_X0identifier}
$\chi$ is an isomorphism. In particular, $X_0$ is smooth, and moreover is isomorphic
to ${\C}^{\dim(G)}$.
\end{thm}

Note that the second claim follows immediately from the first claim
and the well-known fact that a unipotent algebraic group is 
isomorphic to its Lie algebra. 
We will prove several lemmas in order to prove the first claim
in the Theorem. 

The proof we give will essentially follow  \cite{BK},
and is due to DeConcini and Springer \cite{DCS1}.

\begin{lem}
\label{lem_betamap}
There exists a $U^- \times U$-equivariant morphism
$\beta: X_0 \to U^- \times U$ such that 
$\beta \circ \chi (u,v,z)=(u,v)$ for all $z\in \psi(T)$.
\end{lem}

We assume this Lemma for now, and show how to use it to prove
the Theorem.
Recall the following easy fact.

\begin{lem}
\label{lem_principalbundle}
Let $A$ be an algebraic group and let $Z$ be a $A$-variety, and 
regard $A$ as a $A$-variety using left multiplication.
Suppose there exists a $A$-equivariant morphism
$f:Z\to A$. Let $F=f^{-1}(e)$. Then $Z\cong A \times F$.
\end{lem}

\noindent
{\bf Proof.} Define $\chi:A \times F$ by $\chi(a,z)=a\cdot z$,
define $\eta:Z \to F$ by $\eta(z)=f(z)^{-1}\cdot z$, and 
define $\tau:Z\to A\times F$ by $\tau(z)=(f(z),\eta(z))$.
Check that $\chi \circ \tau$ and $\tau \circ \chi$ are
the identity.
\qed

Note also that there is an isomorphism $\sigma: G\times G
\to G\times G$ given by $\sigma (x,y)=(y,x)$. Consider
the isomorphism $\gamma : V\otimes V^* \to V^*\otimes V$ such
that $\gamma (v\otimes f)=f\otimes v$. Then $\gamma\circ a = \sigma(a) \gamma$,
for $a\in G\times G$.

\noindent
{\bf Proof of Theorem \ref{thm_X0identifier}}.
The last Lemma applied to Lemma \ref{lem_betamap} implies
that

\begin{equation}
\label{eq:x0_iso} 
	X_0 \cong (U^- \times U)\times \beta^{-1}(e,e) .
\end{equation}

But $\beta^{-1}(e,e)$ is closed in $X_0$ and by
Lemma \ref{lem_betamap}, $\psi(T)\subset \beta^{-1}(e,e)$.
Thus, $Z=\overline{\psi(T)} \subset \beta^{-1}(e,e)$. 

But $X_0$ is irreducible since $X$ is irreducible, and
$\dim(X_0)=\dim(U^- \times U) + l$, so by \eqref{eq:x0_iso},
$\beta^{-1}(e,e)$ is irreducible and $l$-dimensional.
Thus, $\beta^{-1}(e,e)=Z$, and the proof of the Theorem follows.
\qed

\noindent
{\bf Proof of \ref{lem_betamap}.}
Define a morphism $F:\openp \to \PVopen$ by $F([A])=[A\cdot v_0]$.
The intrinsic definitions of $\openp$ and $\PVopen$ imply that
$F$ is well-defined.

By restriction, we have $\phi:X_0 \to \PVopen$. By Lemma \ref{lem_X0G},
$\psi(U^-TU) \subset X_0$. Let $u\in U^-, t\in T,$ and $v\in U$.
Then $\phi(\psi(utv))=utv \cdot [v_0]=u\cdot [v_0]$ because $TU$
fixes $[v_0]$. Thus, $\phi(\psi(U^-TU))=U^-\cdot [v_0]$, which
is closed in $\PVopen$ by Lemma \ref{lem_uminusclosed}. Since
$\psi(U^-TU)$ is dense in $X_0$, it follows that $\phi(X_0)
= U^-\cdot [v_0]$. It is easy to check that $\phi$ is
$U^-$-equivariant.

Define $\beta_1:X_0 \to U^-$ by setting $\beta_1(x)=u$ if
$u\cdot [v_0]=\phi(x)$. In other words, $\beta_1 = \eta \circ \phi$,
where $\eta$ is the inverse of the action isomorphism $U^- \to U^- \cdot [v_0]$
from Lemma \ref{lem_uminusclosed}. It follows from definitions that
$\beta_1$ is $U^-$-equivariant.

We now define a morphism $\beta_2:X_0 \to U$. To do this,
consider the identifications 
$$
\End(V) \cong V\otimes V^* \cong V^* \otimes V \cong \End(V^*)
$$
where the second morphism is given by $\gamma$. By replacing
$B$ with $B^-$, we may regard $v_0^*$ as a highest weight
vector for $V^*$.
Now the above result applied to $V^*$ and $B^-$ in place of $V$
and $B$
gives a morphism $J:X_0 \to \PVstaropen$ given by $J([A])=
A\cdot [v_0^*]$. Then reasoning as above, $J$ restricts
to a morphism $\psi:X_0 \to U\cdot [v_0^*]$ and then we
set $\beta_2(x)=v\in U$ where $v\cdot [v_0^*] = \psi(x)$.
As before, $\psi$ is $U$-equivariant.

Now $\beta:X_0 \to U^- \times U$ given by $\beta(x)=(\beta_1(x),
\beta_2(x))$ is the required morphism.
\qed 

\subsection{Global Geometry}
\label{sec_globalgeometry}

We use the smoothness of the open set $X_0$ to show $X$ is smooth.

\begin{lem}
\label{lem_uniqueclosedorbit}
Let $K$ be a semisimple algebraic group and let $V$ be an 
irreducible representation of $K$. Then the induced action
of $K$ on $\PR (V)$ has a unique closed orbit through
the line generated by a highest weight vector.
\end{lem}

\noindent
{\bf Proof.} We use the following standard facts:

(1) If
$H$ is a subgroup of $K$, then $K/H$ is projective if and
only if $H$ is parabolic.

(2) A nonzero vector $v$ in
$V$ is a highest weight vector if and only if its stabilizer $K_{[v]}$
in $K$ is parabolic.
If $v_1, v_2$ are highest weight vectors in $V$, then
there is $g\in K$ such that $g\cdot v_1 = v_2$.

\noindent By (1) and (2), it is clear that if $v$ is a highest
weight vector, $K\cdot [v] \cong K/K_{[v]}$ is projective,
and thus closed in $\PR (V)$. Let $[u] \in \PR (V)$ be
such that $K\cdot [u]\cong K/K_{[u]}$ is closed and therefore
projective. Then $K_{[u]}$ is parabolic by (1), so $u$ is
a highest weight vector. Thus, by (2), $K\cdot [u] = K\cdot [v]$.
\qed

\begin{lem}
\label{lem_opencover}
Let $A$ be an algebraic group and let $W$ be a $A$-variety with a
unique closed orbit $Y$. Let $V\subset W$  be open and suppose
$V\cap Y$ is nonempty. Then $W = \cup_{g \in A} \ g\cdot V$.
\end{lem}

\noindent
{\bf Proof.} The union $\cup_{g \in A} \ g\cdot V$ is clearly open,
so $W - \cup_{g \in A} g\cdot V$ is closed, and is $A$-stable.
We assume it is nonempty and argue by contradiction.
Recall that if an algebraic group
acts on a variety, then it has a closed orbit. Thus,
there exists a closed orbit $Z$ for $A$ on
$W - \cup_{g \in A} \  g\cdot V$. Since $Z$ is a closed orbit
for $A$ in $W$, then $Z=Y$. But $Y\not\subset 
W - \cup_{g \in A} \ g\cdot V$ by assumption. It follows that
$W - \cup_{g \in A} g\cdot V$ is empty.
\qed

\begin{lem}
\label{lem_opencoverX}
 Let $W \subset X=\overline{\psi(G)}$ be
a $G\times G$-stable closed subvariety of $X$. Then
$W = \cup_{a\in G\times G} \ a\cdot (W\cap X_0)$.
\end{lem}

\noindent
{\bf Proof.} Note that $\End(V)\cong V\otimes V^*$ is
an irreducible representation of the semisimple group
$G\times G$ and $v_0 \otimes v_0^*$ is
a highest weight vector stabilized by $B\times B^-$.
 Thus, by Lemma \ref{lem_uniqueclosedorbit},
$\PR (\End(V))$ has a unique closed $G\times G$-orbit
through $p_\lambda := [v_0 \otimes v_0^*]$. 
Since $W \subset \PR (\End(V))$ is $G\times G$-stable, it has a
closed orbit, which must be projective since $W$ is projective.
This closed orbit must be $Y=(G\times G) \cdot p_\lambda$. Since
$p_\lambda \subset X_0 \cap Y$, $X_0 \cap W$ is nonempty.
Thus, by Lemma \ref{lem_opencover}, the conclusion
follows.
\qed

\begin{prop}
\label{prop_smoothness}
The following hold:
\begin{enumerate}
	\item $X =  \cup_{a \in G\times G} \ a\cdot X_0$ is smooth.

	\item Let $Q \subset X$ be a $G\times G$-orbit. Then 
$$\overline{Q} = \cup_{a\in G\times G} \ a\cdot (\overline{Q}
\cap X_0).$$ 

	\item If $Q$ and $Q^\prime$ are two $G\times G$-orbits in $X$,
and $\overline{Q} \cap X_0 = \overline{Q^\prime} \cap X_0$,
then $Q = Q^\prime$.
\end{enumerate}

\end{prop}

\noindent
{\bf Proof.} For (2), apply Lemma \ref{lem_opencoverX} to the case $W=
\overline{Q}$. This gives (2), and
in the case $Q=\psi(G)$, we obtain $\overline{Q}=X$, which gives
 $$X= \cup_{a \in G\times G} \ a\cdot X_0 .$$
Now (1) follows by Theorem \ref{thm_X0identifier} since $X$ is a union
of smooth open sets.
To prove (3), note that by (2),
$\overline{Q} = \overline{Q^\prime}$. Since $Q$ and $Q^\prime$
are both open in their closures, they coincide.
\qed

\subsection{Description of the $G\times G$-orbits}
\label{sec_orbitdescription}

We classify the $G\times G$-orbits in $X$ and show they have
smooth closure.

For $I \subset \{ 1, \dots, l \},$ let
$$ Z_I := \{ (z_1, \dots, z_l) \in Z : z_i = 0 \  \forall i\in I \} $$
and let
$$ Z_I^0 := \{ (z_1, \dots, z_l) \in Z_I : z_j \not= 0 \ \forall
j\not\in I \}.$$
Then $Z_I$ is the closure of $Z_I^0$, 
$Z_I \cong {\C}^{l - |I|}$ and $Z_I^0 \cong {(\cstar)}^{l - |I|}$.
Moreover, by the proof of Lemma \ref{lem_Zaffine}, 
$T\cong {(\cstar)}^l$ acts on $Z\cong \C^l$ by $(a_1, \dots, a_l)\cdot (z_1, \dots, z_l)
=(a_1z_1, \dots, a_lz_l)$ in appropriate coordinates, so 
$Z$ is $(T \times \{ e \})$-stable, and the $Z_I^0$
are exactly the $T \times \{ e \}$-orbits in $Z$.

The next lemma follows from the above discussion.

\begin{lem}
\label{lem_zIorbit}
Set $z_I = (z_1, \dots, z_l )$, where $z_i = 1 \ {\rm if \ } i \not\in I$, and
$z_i = 0 \ {\rm if \ } i \in I$.
Then $ Z_I^0 = (T \times \{ e \})\cdot z_I$.
\end{lem}

For $i = 1, \dots, l$, $Z_i := Z_{\{ i \}}$ is a hypersurface,
and $Z_I = \cap_{i\in I} Z_i$.

\noindent Let 
$$\Sigma_I = \chi(U^- \times U \times Z_I) $$
and let
$$\Sigma_I^0 = \chi(U^- \times U \times Z_I^0)$$

The following result follows easily from Theorem \ref{thm_X0identifier}.

\begin{prop}
\label{prop_xoorbits}
The $\Sigma_I^0$ are precisely the $U^-T \times U$ orbits
in $X_0$, and the closure $\Sigma_I$ of $\Sigma_I^0$ in
$X_0$ is isomorphic to ${\C}^{\dim(G) - |I|}$. In particular,
$U^-T \times U$ has precisely $2^l$ orbits on $X_0$, and all these
orbits have smooth closure.
\end{prop}

\begin{rem}
\label{rem_x0vsxorbits}
Let $Q$ be a $G\times G$-orbit in $X$. Then
 $\overline{Q} \cap X_0$ is closed, irreducible,
and $U^-T \times U$-stable, so it must be the closure of one of the $2^l$
$U^-T \times U$-orbits in $X_0$.
\end{rem}

By this last remark
 and Proposition \ref{prop_smoothness}
(3), it follows that $G\times G$ has at most $2^l$ orbits on $X$.
We show there are exactly $2^l$ $G\times G$ orbits.

For the following result, see \cite{hartshorne}, chapter II.

\begin{lem}
\label{lem_affinecomplement}
Let $W$ be an irreducible projective variety with open set $U$. Then if
$U$ is isomorphic to an affine variety, all irreducible
components of $W - U$ have codimension $1$.
\end{lem}

\begin{lem}
\label{lem_gcomplement}
$X - \psi(G) = \cup_{i = 1, \dots, l} \ S_i$, where 
\begin{enumerate}
	\item $S_i$ is a $G\times G$-stable divisor 

	\item $S_i = \overline{\Sigma_i}$, and

	\item $S_i \cap X_0 = \Sigma_i$.
\end{enumerate}
\end{lem}

\noindent
{\bf Proof.}
Let $X - \psi(G) = \cup_{\alpha} \ S_\alpha$ be the decomposition of
$X - \psi(G)$ into irreducible components.
Since $\psi(G)\cong G$ is affine, 
each $S_\alpha$ has codimension one by Lemma \ref{lem_affinecomplement}.
 Since $X - \psi(G)$ is $G\times G$-stable,
and $G\times G$ is connected, each $S_\alpha$ is $G\times G$
stable. Thus, $S_\alpha \cap X_0$ is $U^-T \times U$-stable
and is a closed irreducible hypersurface in $X_0$ by Lemma
\ref{lem_opencoverX}. Since $X_0$ has
finitely many $U^-T \times U$-orbits, it must be one of the
$\Sigma_i$. Conversely, the closure of $\Sigma_i$ in $X$
is a closed, irreducible hypersurface. Since $X_0 \cap
\psi(G) = \Sigma_{\emptyset}$ by Lemma \ref{lem_X0G},
$\Sigma_i$ is contained in the closed set $X - \psi(G)$,
so $\overline{\Sigma_i}$ is contained in $X - \psi(G)$.
Thus, $\overline{\Sigma_i}$ must be an irreducible component
of $X - \psi(G)$, so $\overline{\Sigma_i} = S_\alpha$ for
some $\alpha$, and $S_\alpha \cap X_0 = \Sigma_i$. The result
follows.
\qed

\begin{lem}
\label{lem_SIdef}
The following hold:

\begin{enumerate}
	\item Let $S_I = \cap_{i \in I} S_i$. Then $S_J \subset S_I$ 
	if $I \subset J$ and $S_I \cap X_0 = \Sigma_I$.

	\item Let $S_I^0 = S_I - \cup_{I \subsetneqq J} \ S_J$. Then $S_I^0 = 
	(G\times G)\cdot z_I$ is a single $G\times G$-orbit and $S_I^0 = S_J^0$
	implies $I=J$.

	\item $S_I = \cup_{a\in G\times G}\  a\cdot \Sigma_I$, and in 
	particular $S_I$ is smooth.
\end{enumerate}
\end{lem}

\noindent
{\bf Proof.} The first claim of (1) is obvious from the
definition and the second claim follows from Lemma \ref{lem_gcomplement}(3).
 For (2), it is
clear that $S_I^0 \cap X_0 = \Sigma_I^0$. Since the
$S_i$ are $G\times G$-stable, $S_I$ is $G\times G$-stable,
so $S_I^0$ is $G\times G$-stable. Let $x, y
\in S_I^0$. By Proposition \ref{prop_smoothness} (1), there are $a$ and
$b$ in $G\times G$ such that $a\cdot x$ and $b\cdot y$
are in $X_0$. Since $a\cdot x$ and $b\cdot y$ are in
$S_I^0$ by $G\times G$-stability, they are in $\Sigma_I^0$.
Thus, by Proposition \ref{prop_xoorbits}, there is
$c\in U^-T \times U$ such that $a\cdot x = cb \cdot y$.
Part (2) follows since clearly $z_I \in S_I^0$.
Part  (3) follows by the above assertions
and Proposition \ref{prop_smoothness}(2). 
\qed

\begin{rem}
\label{rem_normalcrossings}
Recall the following definition. Let $X$ be a smooth variety
with hypersurface $Z$. We say $Z$ is a divisor with normal
crossings at $x\in Z$ if there is an open neighborhood $U$
of $x$ such that $Z\cap U = D_1 \cup \dots \cup D_k$ is
a union of hypersurfaces and $D_{i_1} \cap \dots \cap
D_{i_j}$ is smooth of codimension $j$ for each distinct $j$-tuple
$\{ i_1, \dots, i_j \}$ in $\{ 1, \dots, k \}$. A standard example is
to take $X = {\C}^n$ and $Z$ to be the variety defined by
the vanishing of $z_1 \dots z_k$. Then $Z$ is the union of
the hyperplanes given by vanishing of $z_i$, and in the complex
analytic setting, every
divisor with normal crossings is locally of this nature.
\end{rem}

\begin{thm}
\label{thm_orbitclosure}
$G\times G$ has $2^l$ orbits in $X$, given by
 $S_I^0 = (G\times G)\cdot z_I$ where
$I$ is a subset of $\{ 1, \dots, l \}$. In particular, all orbits
have smooth closure, and the pair $(X,X-\psi(G))$
is a divisor with normal crossings.
\end{thm}

\noindent
{\bf Proof.}
 Let $Q \subset X$
be a $G\times G$-orbit. Then $\overline{Q} \cap X_0 = 
\Sigma_I$ for some $I \subset \{ 1, \dots, l \}$ by
Remark \ref{rem_x0vsxorbits} and Proposition \ref{prop_xoorbits}.
It follows easily that $\Sigma_I^0 \subset Q$ so $S_I^0=Q$, using
Lemma \ref{lem_SIdef} (2). Again by Lemma \ref{lem_SIdef} (2), there
are exactly $2^l$ orbits. 
Since $\Sigma_I$ is the closure of 
$\Sigma_I^0$ in $X_0$, it follows that $S_I$ is the closure
of $S_I^0$ in $X$ using Lemma \ref{lem_opencover}. 
It is clear from Propositon \ref{prop_xoorbits}
 that $(X_0,X_0 - \psi(U^-TU))$
is a divisor with normal crossings. The last assertion follows
from Proposition \ref{prop_smoothness} (1) 
and group action invariance of the divisor with normal
crossings property.
\qed

\subsection{Geometry of orbits and their closures}
\label{sec_geometryoforbits}

We now want to understand the geometry of $S_I^0$ and
$S_I$. We show that the orbit $S_I^0$ fibers over a product of generalized
flag varieties, with fiber a semisimple group. In this picture, the closure $S_I$
 fibers over the same product of generalized
flag varieties, and the fiber is the wonderful compactification of
the semisimple group.

For a subset $I \subset \{ 1, \dots, l \}$, let
$\Delta_I = \{ \alpha_i : i \not\in I \}$. Let
$\Phi_I$ be the roots that are in the linear
span of the simple roots in $\Delta_I$, and let
$$
\fl_I = \ft + \sum_{\alpha \in \Phi_I} \fg_\alpha
$$
Let $\fp_I = \fl_I + \fu$ and let $\fp_I^- = \fl_I + \fu^-$.
Let $\fu_I$ and $\fu_I^-$ be the nilradicals of
$\fp_I$ and $\fp_I^-$, respectively. Note that
$\fp_\emptyset = \fg$ and $\fp_{\{1, \dots, l \}} = \fb$.
Let $P_I$, $P_I^-$, $U_I$, $U_I^-$, and $L_I$ be the
 corresponding connected subgroups of $G$. Let $Z(L_I)$
be the center of $L_I$ and let $G_I = L_I/Z(L_I)$ be the
adjoint group of $L_I$ and denote its Lie algebra by
$\fg_I$.

% Let $\tilde{P_I), \tilde{L_I}, \tilde{Z(L_I)}$, etc.,
%be the preimages of these groups in $\tilde{G}$.

For a Lie algebra $\fa$, let $U(\fa)$ be its enveloping
algebra. Let $V_I = U(\fl_I)\cdot v_0$. Then $V_I$ is a
$\fl_I$-stable submodule of $V$ and $V_I$ may be 
identified with the irreducible representation of
$\fg_I$ of highest weight $\lambda$. 

\begin{lem}
\label{lem_VIaction} The following hold:
\begin{enumerate}
	\item $V_I$ is $\fp_I$-stable, and $\fu_I$ annihilates $V_I$. 

	\item Let $Q = \{ g \in G : g\cdot V_I = V_I$\}. Then $Q=P_I$. 
\end{enumerate}
%
%let $V^{\fu_I^-}$ be the subspace
%of $V$ annihilated by $\fu_I^-$. Then
%$V^{\fu_I^-} \cap V_I = \{ 0 \}$.
\end{lem}

\noindent
{\bf Proof.} For (1), by the Poincare-Birkhoff-Witt Theorem,
$U(\fp_I)=U(\fl_I)U(\fu_I)$, and $\fp_I$-stability
of $V_I$ follows by the theorem of the highest weight.
Since $\fu_I$ is an ideal of $\fp_I$, the space
of $\fu_I$-invariants $V^{\fu_I}$ is $\fl_I$-stable and meets $V_I$.
Since it is nonzero and $V_I$ is an irreducible 
$\fl_I$-module, $V_I = V^{\fu_I}$.

For (2), note that it follows from (1) that $P_I \subset Q$. Hence
 $Q$ is parabolic so it is
connected. To show $P_I = Q$, it suffices to show
$\fp_I = \fq$, where $\fq$ is the Lie algebra of
$Q$. There is $t\in Z(L_I)$  such that $\alpha_i(t)\not= 1$
 for every simple root $\alpha_i \not\in \Delta_i$.
 Then if $\tilde{t}\in \tilde{T}$ is a preimage of $t$, ${Ad}_{\tilde{t}}$
acts by $\lambda(\tilde{t})$ on $V_I$, and if $x\in \fg_{-\alpha_i}$
for $\alpha_i$ simple, ${Ad}_{\tilde{t}}$ has
weight $\frac{\lambda(\tilde{t})}{\alpha_i(t)}$ on 
$x\cdot V_I$. 
It follows that if $\alpha_i \not\in \Delta_i$,
 $\fg_{-\alpha_i}\cdot  V_I \cap V_I = 0$.
Since $\lambda$ is regular, $\fg_{-\alpha_i}\cdot v_0 \not= 0$ for each
simple root $\alpha_i$, so $\fg_{-\alpha_i}\cdot v_0 \not\in V_I$ if
$\alpha_i \not\in \Delta_I$. It follows
that $\fg_{-\alpha_i} \not\subset \fq$ for all $\alpha_i \not\in  \Delta_I$, so
$\fq=\fp_I$.
\qed

\pagebreak

\begin{lem}
\label{lem_VIbasis} The following hold:
\begin{enumerate}
	\item Let 
	$$J = \{ k\in \{ 0, \dots, n \} : \lambda_k = \lambda - 
	\sum_{ i \in \Delta_I } n_i \alpha_i, n_i \ge 0 \}
	$$
	Then the $\{ v_j : j\in J \}$ form a basis of $V_I$.

	\item Recall $z_I = \{ (\epsilon_1, \dots, \epsilon_l)
	\in Z_I \}$, where $\epsilon_i = 1$ if $i \not\in I$
	and $\epsilon_i = 0$ if $i\in I$. Then $z_I
	= [{\pr}_{V_I}]$, where ${\pr}_{V_I}$ is the projection on
	$V_I$ such that ${\pr}_{V_I}(v_k)=0$ if $k\not\in J$.
\end{enumerate}
\end{lem}

\noindent
{\bf Proof.} (1) is routine and is left to the reader.
For (2), recall that for $\{ z_1, \dots, z_l \} \in Z$,
the corresponding class in $\PR (\End(V))$ is 
$$[v_0 \otimes v_0^* + \sum_{i=1}^l z_i v_i \otimes v_i^*
+ \sum_{k > l} \prod_{i=1}^l z_i^{n_{ik}} v_k\otimes v_k^*] $$
where $\lambda_k = \lambda - \sum_{i=1}^l n_{ik} \alpha_i$.
It follows that
$$z_I = [ \sum_{k \in J}  v_k \otimes v_k^*].$$
%where $c_k = 1$ if $k \in J$ and $c_k = 0$ if $k \not\in J$.
But $\sum_{k \in J}  v_k \otimes v_k^*$ is easily identified
with ${\pr}_{V_I}$.
\qed

We now compute the stabilizer of $(G\times G)$ at the point
$z_I \in S_I^0 = (G\times G)\cdot z_I$.

\begin{prop}
\label{prop_zIstabilizer}
The stabilizer $(G\times G)_{z_I}$ of $z_I$ in $G\times G$ is 
$$\{ (xu, yv) : u \in U_I, v\in U_I^-, x\in L_I, y\in L_I,
\ {\rm and } \ x y^{-1} \in Z(L_I) \}$$
In particular, $S_I^0$ fibers over $G/P_I \times G/P_I^-$
with fiber $G_I$.
\end{prop}

\noindent
{\bf Proof.} Suppose $(r,s)\cdot z_I=z_I$ for $(r,s)\in G\times G$.
 Then
$[r{\pr}_{V_I} s^{-1}]=[{\pr}_{V_I}]$ so $r$ preserves the image
$V_I$ of ${\pr}_{V_I}$. Thus, by Lemma \ref{lem_VIaction},
$r\in P_I$. Further, if $r\in U_I$, then $r\cdot [{\pr}_{V_I}]
= {\pr}_{V_I}$ since $U_I$ acts trivially on $V_I$.

Recall the identification $\gamma :\End(V)\to  
\End(V^*)$ from the proof of \ref{lem_betamap}.
 Let $V_I^* = U(\fp^-)v_0^* \subset V^*$. $V_I^*$ is an
irreducible representation of $L_I$ with lowest weight $-\lambda$.
One can easily check that $\gamma({\pr}_{V_I})
={\pr}_{V_I^*}$, and
 $\gamma ((x,y)\cdot A) = (y,x)\cdot \gamma(A)$, for
$(x,y)\in G\times G$ and $A\in \End(V)$. Hence
$$
(s,r)\cdot [{\pr}_{V_I^*}] = \gamma((r,s) \cdot [{\pr}_{V_I}])=
\gamma([{\pr}_{V_I}])=[{\pr}_{V_I^*}].
$$
As above, it follows that
 $s$ preserves $V_I^*$. Thus, by Lemma \ref{lem_VIaction}
applied to the opposite parabolic, $s \in P_I^-$. Further, if
$s\in U_I^-$, then $[{\pr}_{V_I} \circ s]=[{\pr}_{V_I}]$.

Now let $r=xu$ and let $s=yv$ with $x, y \in L_I$, $u\in U_I$
and $v\in U_I^-$. Then 
$$r\cdot [{\pr}_{V_I}] \cdot s^{-1}= x \cdot [{\pr}_{V_I}] \cdot y^{-1}.$$
In particular, $xy^{-1}$ acts trivially on $\PR (V_I)$. But
$$
Z(L_I)= \{ g \in L_I : g\cdot [v]=[v], \forall [v] \in \PR (V_I) \},
$$
so $xy^{-1}\in Z(L_I)$.
Now we consider the projection $(G \times G)/(G\times G)_{z_I}
\to (G\times G)/(P_I \times P_I^-)$ with fiber
$$(P_I \times P_I^-)/(G\times G)_{z_I} \cong 
(L_I \times L_I)/\{ (x, y) :  x\in L_I, y\in L_I,
\ {\rm and } \ x y^{-1}\in Z(L_I) \}
$$
The morphism $$(L_I \times L_I)/\{ (x, y) :  x\in L_I, y\in L_I,
\ {\rm and } \ x y^{-1}\in Z(L_I) \} \to L_I/Z(L_I)$$
given by $(a,b)\mapsto ab^{-1}$ is easily seen to be an
isomorphism.
\qed

To understand the closure $S_I$ of $S_I^0$, we embed the compactification
of a smaller group into $X$.

First, note that we may embed $\End(V_I)$ into $\End(V)$
by using the map
$$ 
\sum_{i, j\in J} a_{ij} v_i \otimes v_j^* \to  
\sum_{i, j\in \{ 0, \dots, n \} } b_{ij} v_i \otimes v_j^*$$
where $b_{ij}=a_{ij}$ if $i,j \in J$, and $b_{ij}=0$ otherwise.

This map induces an embedding $\PR (\End(V_I)) \to \PR (\End(V))$,
and we will regard $\PR (\End(V_I))$ as a closed subvariety
of $\PR (\End(V))$.

It follows from definitions that $z_K \in \PR (\End(V_I)) \iff I \subset K$.

Define a morphism $L_I \to \PR (\End(V))$ by $\psi_I (g)=[g {\pr}_{V_I}]$.
Then $\psi_I$ descends to a morphism $\psi_I:L_I/Z(L_I)=G_I
\to \PR (\End(V))$. Note that the image of $\psi_I$ is in 
$\PR (\End(V_I))$. 

Recall that the center of $G_I$ is trivial,
and note that $V_I$ is an irreducible representation of
a cover of $G_I$ with regular highest weight. We set
$X_I = \overline{\psi_I(G_I)}$, and we may apply our results about
orbit structure for the compactification $X$ of $G$ to the compactification
$X_I$ of $G_I$. In particular, $X_I$ is a $G_I \times G_I$-variety
and we may regard $X_I$ as a $P_I \times P_I^-$-variety
via the projection $P_I \times P_I^- \to (P_I \times P_I^-)/
(Z(L_I)U_I \times Z(L_I)U_I^-)=G_I \times G_I$.

\begin{thm}
\label{thm_GIcompactification}
Consider the morphism
$$
\chi:(G\times G)\times_{P_I \times P_I^-} X_I \to S_I,
\ \chi(g_1,g_2,x)=(g_1,g_2)\cdot x.$$
Then, $\chi$ is an isomorphism of varieties. In particular,
$S_I$ fibers over $G/P_I \times G/P_I^-$ with fiber $X_I$.
\end{thm}

\noindent {\bf Proof.}
It suffices to prove $\chi$ is a bijection, since a bijection
to a smooth variety in characteristic zero is an isomorphism.
For $I\subset K$, let 
$y_K=(e,e,z_K) \in (G\times G)\times_{P_I \times P_I^-} X_I$
and note that $\chi(y_K)=z_K$. Since $\chi$ is $G\times G$-equivariant,
it follows that $\chi$ is surjective by Theorem \ref{thm_orbitclosure}. 

To show $\chi$ is injective, we use the formal  fact:

(*) Let $A$ be a group with subgroup $B$. Let $Y$ be
a $B$-set. Then $B$-orbits in $Y$ correspond to $A$-orbits in $A\times_B Y$
via $B\cdot y \mapsto A\cdot (e,y)$. Moreover, the stabilizer
$B_y$ coincides with the stabilizer $A_{(e,y)}$.

It follows that the stabilizer $(G\times G)_{y_K}=
(P_I \times P_I^-)_{z_K}$. 
But $U_I Z(L_I) \times U_I^- Z(L_I)$ acts trivially
on $X_I$ since it  fixes each $z_K$ with $I\subset K$,
so $(P_I \times P_I^-)_{z_K}$ is the preimage
of $(G_I \times G_I)_{z_K}$ in $P_I \times P_I^-$.
By Proposition \ref{prop_zIstabilizer} applied to $G_I$ and $X_I$,
it follows that
$$
(P_I \times P_I^-)_{z_K} = 
\{ (xu, yv) : u \in U_K, v\in U_K^-, x\in L_K, y\in L_K,
\ {\rm and } \ x y^{-1} \in Z(L_K) \},$$
which coincides with the stabilizer $(G\times G)_{z_K}$.
Thus, $\chi$ is injective when restricted to a $G\times G$ orbit.

To complete the proof that $\chi$ is injective, it suffices to show
that every $G\times G$-orbit contains some $y_K$, $I\subset K$.
This is an easy consequence of (*) and the classification of
$G_I$ orbits in its wonderful compactification $X_I$, which follows
from Theorem \ref{thm_orbitclosure}. We leave the routine details
to the reader.
\qed

\section{Comparison of different compactifications}
\label{sec_comparisonofcompactifications}

\subsection{Independence of highest weight}
\label{sec_independence}

For completeness, we show that $X$ is independent of the
choice of the regular, dominant highest weight. This result
is not especially surprising since $\lambda$ does not
appear in the statements describing the $G\times G$-orbit
structure. 
 The proof mostly 
follows that of DeConcini and Springer (\cite{DCS1}, Proposition 3.10).

 We may 
carry out the wonderful compactification construction using regular, dominant weights 
$\lambda_1$ and $\lambda_2$ to get varieties $X^1 \subset \PR( \End V(\lambda_1) )$ and
 $X^2 \subset \PR( \End V(\lambda_2) )$ respectively. 
We denote the class of the identity by $[id_1] \in X^1$ and $[id_2] \in X^2$.
We define $X^\Delta = \overline{ (G \times G) \cdot ([id_1], [id_2] ) } \subset X^1 \times X^2$,
and prove that the natural projections $p_1 : X^\Delta \to X^1$ and $p_2: X^\Delta \to X^2$
are isomorphisms. 

\begin{prop} 
\label{prop_independence}
For $i=1, 2$, $p_i: X^{\Delta} \to X^i$ is a $G \times G$-equivariant isomorphism of varieties. 
In particular, $p_2 \circ p_1^{-1}: X^1 \to X^2$ is a $G \times G$-equivariant isomorphism sending $[id_1]$ to 
$[id_2]$.
\end{prop}

In the proof, we consider constructions in $X^i$ analogous 
 to constructions used in Section \ref{sec_construction}, such as
the open affine piece, and the closure of $T$ in the open affine piece. We denote the
various subsets defined in Section \ref{sec_construction} for $X^i$ with a 
superscript $i$, e.g., $Z^i \subset X^i_0 = \openp( \End V(\lambda_i) ) \cap X^i$. 

First, let $Z^{\Delta} := \overline{ (T \times \{e\})\cdot([id_1], [id_2]) }$ where closure is 
taken in the open set $X^1_0 \times X^2_0 \subset X^1 \times X^2$. 

\begin{lem} \label{lem_Zdelta}
$Z^\Delta \cong \mathbb{C}^l$ and $p_i: Z^\Delta \to Z^i$ is a $(T \times \{e\})$-equivariant isomorphism.
\end{lem}

\noindent {\bf Proof.}
This is a straightforward calculation using coordinates as in Lemma \ref{lem_Zaffine}. 
Choose a basis of weight vectors  $\{v_i\}_{i=0..n}$ for $V(\lambda_1)$
and a basis of weight vectors $\{w_i\}_{i=0..m}$ for $V(\lambda_2)$
satisfying properties (1)-(3) preceding Remark \ref{rem_weightmultiplicity}.
 Then, for $t \in T$,
\[ 
\begin{array}{ll}(t \times \{e\})\cdot([id_1], [id_2]) & = (t \times \{e\})\cdot( [ \sum_{i=0}^n v_i^* \otimes v_i ], 
[\sum_{i=0}^m w_i^* \otimes w_i ]) \\
 & = ( [ v_0^* \otimes v_0 + \sum_{i=1}^l \frac{1}{\alpha_i(t)} v_i^*
 \otimes v_i + \cdots ], \\
& \phantom{ = ( }[ w_0^* \otimes w_0 + \sum_{i=1}^l \frac{1}{\alpha_i(t)} w_i^* \otimes w_i + \cdots ] )
 \end{array}
 \]
 
\noindent where the terms indicated by $\cdots$ in the two factors have coefficients which are polynomial in the 
$\frac{1}{\alpha_i(t)}$. Then the map sending  $(z_1, \dots, z_l) \in \mathbb{C}^l$
 to the above expression 
with $\frac{1}{\alpha_i(t)}$ replaced by $z_i$ identifies $Z^\Delta$ with $\mathbb{C}^l$ and the  
claim follows.
\qed
 
Recall the embedding from Theorem \ref{thm_X0identifier},
$$\chi^i : U^- \times U \times Z^i \to X^i$$
This is 
an isomorphism to $X^i_0$. 
Let $X^\Delta_0 := p_i^{-1}(X^i_0) \subset X^\Delta$, and consider the subset
$V= \bigcup_{a \in G \times G} \ a\cdot X^\Delta_0$ of $X^\Delta$. Define $\chi^\Delta : U^- \times U \times Z^\Delta \to X^\Delta$ in the same manner as $\chi$ from Theorem \ref{thm_X0identifier}.

We have a commutative diagram:
\begin{equation} \label{eq:indep_square}
\xymatrix{
U^- \times U \times Z^\Delta \ar[d]^{id \times id \times p_i}_{\cong} \ar[r]^-{\chi^\Delta} & X^\Delta 
\ar[d]^{p_i} & \\
U^- \times U \times Z^i \ar[r]^-{\chi^i} & X^i & .
} 
\end{equation}
We will use the following two lemmas, which we prove below.
\begin{lem}
\label{indep_lemma_embed}
$\chi^\Delta$ is an embedding onto the open set $X^\Delta_0$.
\end{lem}

\begin{lem}
\label{indep_lemma_inj}
$p_i \lvert_{V}$ is injective.
\end{lem}

\noindent
We assume these for now and prove that $p_i$ is an isomorphism. 

\noindent {\bf Proof  of Proposition \ref{prop_independence}.}
Consider the commutative diagram:
\begin{equation} \label{eq:indep_square_restricted}
\xymatrix{
U^- \times U \times Z^\Delta \ar[d]^{id \times id \times p_i}_{\cong} \ar[r]^-{\chi^\Delta} & X_0^\Delta 
\ar[d]^{p_i} & \\
U^- \times U \times Z^i \ar[r]^-{\chi^i} & X_0^i & .
}
\end{equation}

\noindent  Since $\chi^i:U^- \times U \times Z^i \to X_0^i$
is surjective, $p_i:X_0^\Delta \to X_0^i$ is surjective
using the above commutative diagram. Then $p_i:V\to X^i$ is surjective
by Proposition \ref{prop_smoothness} (1). Further, $p_i:V \to X^i$
is injective by Lemma \ref{indep_lemma_inj}, and hence an isomorphism
since $X^i$ is smooth. 

\pagebreak

Thus $V$ is
complete, so $V=X^\Delta$, since $X^\Delta$ is irreducible
 and $\dim(V)=\dim(X^\Delta)$.\qed

We now prove the two lemmas.

\noindent {\bf Proof of Lemma \ref{indep_lemma_embed}.}
First, $\chi^\Delta$ is an embedding since
 $\chi^i \circ (id \times id \times p_i\lvert_{Z^\Delta})$ is
an embedding and the diagram \eqref{eq:indep_square} commutes.

Let $Y$ denote the image of $\chi^\Delta$. Commutativity of \eqref{eq:indep_square} and Theorem \ref{thm_X0identifier} imply that $Y \subset X_0^\Delta$. Consider the map 
$$s = \chi^\Delta \circ (id \times id \times p_i)^{-1} \circ (\chi^i)^{-1} : X_0^i \to X_0^\Delta .$$
Note that $s$ is a section for the map $p_i$ over $X_0^i$, i.e., $p_i \circ s = id$. 
Consider the composition $f = s \circ p_i\lvert_{X_0^\Delta}$. It is routine
to check that 
$f\lvert_Y = id\lvert_Y$. 
 $Y$ and $X_0^\Delta$ have the same dimension and $X_0^\Delta$ is
irreducible since it is an open subset of the irreducible variety $X^\Delta$.
Hence, there is an open, dense subset of $X_0^\Delta$ contained in $Y$. Then
$f$ is the identity on an open, dense subset of $X_0^\Delta$, so $f$ is the identity
on $X_0^\Delta$. Hence, $Y=X_0^\Delta$.
\qed

\noindent {\bf Proof of Lemma \ref{indep_lemma_inj}.}
For $I \subset \{1, \ldots, l\}$, let $z^\Delta_I = (z_I^1, z_I^2)$,
so $p_i(z^\Delta_I)=z_I^i$. Then the stabilizer
$(G \times G)_{z^\Delta_I} = (G \times G)_{z^1_I}\cap (G \times G)_{z^2_I}$,
which coincides with $(G \times G)_{z^i_I}$ for $i=1, 2$,
since $(G \times G)_{z^1_I} = (G \times G)_{z^2_I}$ by Proposition
\ref{prop_zIstabilizer}. Thus, $p_i:V \to X^i$ is injective when restricted
to the orbits through some $z^\Delta_I$. By Lemma \ref{indep_lemma_embed}, each
$U^-T \times U$-orbit on $X_0^\Delta$ meets $Z^\Delta$, so it meets some
$z_I^\Delta$. The Lemma follows.
\qed

We now prove a result related to Proposition \ref{prop_independence} which will be used in the sequel.

Let $E$ be a representation of $\tilde{G}\times \tilde{G}$.
Then the $\tilde{G}\times \tilde{G}$ action on $\PR (E)$ descends
to an action of $G\times G$ on $\PR (E)$. Suppose there exists
a point $[x]\in \PR(E)$ such that $(G\times G)_{[x]}=G_\Delta$.
We may embed $G$ into $\PR(E)$ by the mapping
$\psi:G \to \PR(E)$ given by $\psi(g)=(g,e)\cdot [x]$.
Let $X(E,[x])$ be the closure of $\psi(G)$ in
$\PR(E)$. Let $X_{\lambda}=X(\End(V),[id_V])$ when $V$ is
irreducible of
highest weight $\lambda$. If $\lambda$ is regular, then
$X_\lambda$ is of course smooth and projective with known
$G \times G$-orbit structure by Theorem
\ref{thm_orbitclosure}.

Let $W_1, \dots, W_k$ be a collection of irreducible representations
of $G$ of highest weights $\mu_1, \dots, \mu_k$. Let 
$W = W_1 \oplus \cdots \oplus W_k$.
Let $F$ be a representation of $\tilde{G} \times \tilde{G}$. Let
$V$ have highest weight $\lambda$ as before. When useful, we
will denote the irreducible representation of highest weight
$\lambda$ by $V(\lambda)$.

\begin{prop}
\label{prop_lowerweightcomp}
Suppose each $\mu_j$ is of the form 
$\mu_j = \lambda - \sum n_i \alpha_i$ with all $n_i$ nonnegative
integers. Then,

$X(\End(V)\oplus \End(W) \oplus F, [{\id}_V \oplus {\id}_W \oplus 0])
\cong X_\lambda$.
\end{prop}

\noindent {\bf Proof.}
It is immediate from definitions
that 
$$X(\End(V)\oplus \End(W) \oplus F, [{\id}_V \oplus {\id}_W \oplus 0])
\cong X(\End(V)\oplus \End(W), [{\id}_V \oplus {\id}_W ]),$$
i.e., nothing is lost by setting $F=0$. Indeed, the $G\times G$-orbit
through $[{\id}_V \oplus {\id}_W \oplus 0]$ lies inside the
closed subvariety $\PR (\End(V)\oplus \End(W))$ of
$\PR (\End(V)\oplus \End(W) \oplus F)$, and this implies the claim.

Define $X^\prime =  X(\End(V)\oplus \End(W), [{\id}_V \oplus {\id}_W ])$
for notational simplicity.

Consider the open subset 
$\tildePR (\End(V)\oplus \End(W))=\{ [A\oplus B]: A\not= 0 \}$
of $\PR (\End(V)\oplus \End(W))$. The projection
$\pi:\tildePR (\End(V)\oplus \End(W)) \to \PR (\End(V))$
given by $[A\oplus B] \mapsto [A]$ is a morphism of varieties.

We claim  that $\tildePR (\End(V)\oplus \End(W))$ is $(G\times G)$-stable.
Indeed, if $A\not= 0$ and $(x,y)\in \tilde{G}\times \tilde{G}$,
then $xAy^{-1}$ is nonzero since
$x$ and $y$ are invertible. It follows that if
$[A\oplus B] \in \tildePR (\End(V)\oplus \End(W)),$ then
$(x,y)\cdot [A\oplus B] \in \tildePR (\End(V)\oplus \End(W))$,
which establishes the claim.

Let $w_0(j)$ be a highest weight vector for $W_j$, and
let $w_0(j)^*$ be a nonzero vector of weight $-\mu_j$ 
of $W_j^*$, normalized so $w_0(j)^*(w_0(j))=1$.
Let $w_0 = \sum_{j=1}^k w_0(j)$ and
let $w_0^* = \sum_{j=1}^k w_0(j)^*$. Then $w_0^*(w_0)=k$.
Define an
open subset $\PR_0(\End(V)\oplus \End(W))$ as the set
of $[A\oplus B]$ such that $v_0^*(A\cdot v_0)\not= 0$
and $w_0^*(B\cdot w_0)\not= 0$. It is easy to see
that $\PR_0(\End(V)\oplus \End(W))$ is $T\times T$-stable.

Let $Z^\prime$ be the closure of $\psi(T)$ 
in $\PR_0(\End(V)\oplus \End(W))$.

Claim: $Z^\prime \cong \C^l $.

The proof of this claim is essentially the same as the proof
of Lemma \ref{lem_Zaffine}. Indeed, we may compute
$\psi(t)=(t,e)\cdot [{\id}_V \oplus {\id}_W]$ in the same
manner as in Lemma \ref{lem_Zaffine}. In \ref{lem_Zaffine},
the coordinates $z_i$ are essentially given by
$\frac{1}{\alpha_i(t)}$ and the assumption that $\mu_j \le \lambda$ implies
that the additional
summands that appear in $(t,e)\cdot [{\id}_V \oplus {\id}_W]$
have coefficients that are polynomial in the $z_i$. We leave details
to the reader.

It is easy to see that $Z^\prime \subset 
\tildePR (\End(V)\oplus \End(W))$, and in fact $\pi:Z^\prime
\to Z$ is an isomorphism compatible with the identifications
with $\C^l$. It follows from $(G\times G)$-stability
of $\tildePR (\End(V)\oplus \End(W))$ that $X_1^\prime :=
(G\times G)\cdot Z^\prime$ is in 
$\tildePR (\End(V)\oplus \End(W))$.

By Theorem \ref{thm_orbitclosure},  $X = (G\times G)\cdot Z$. It follows that
$\pi:X_1^\prime \to X$ is surjective.

We claim that $\pi:X_1^\prime \to X$ is injective, so
that $X^\prime \cong X$ since $X$ is smooth.
It follows that $X_1^\prime$ is projective since $X$ is projective,
so $X_1^\prime$ is closed in $X^\prime$, and has dimension
equal to the dimension of $G$. But $X^\prime$
is irreducible of dimension equal to the dimension of $G$
so $X^\prime = X_1^\prime$.

For $I \subset \{ 1, \dots, l \}$, define $z_I^\prime$
by the same formula as in Lemma \ref{lem_zIorbit}.
It is routine to show that $Z^\prime =\cup_I \ (T,e)\cdot z_I^\prime$
and $\pi(z_I^\prime)=z_I$.
To show that $\pi$ is injective, it suffices to check that the 
stabilizer $(G\times G)_{z_I^\prime}
=(G\times G)_{z_I}.$ From $(G\times G)$-equivariance of
$\pi$, it follows that 
$(G\times G)_{z_I^\prime} \subset (G\times G)_{z_I}.$
We computed $(G\times G)_{z_I}$ in Proposition \ref{prop_zIstabilizer}.
Using this computation, it is not difficult to check that
$(G\times G)_{z_I} \subset (G\times G)_{z_I^\prime}.$
Indeed, we have
$ z_I^\prime = [ {\pr}_{V_I} + \sum_{j=1}^k {\pr}_{W_{j,I}} ]$,
where $W_{j,I}$ is $U(\fl_I)\cdot w_0(j)$.
By Proposition \ref{prop_zIstabilizer}, the stabilizer $(G\times G)_{z_I}$
acts by scalars on $U(\fl_I)\cdot w_0(j)$. 
This implies injectivity of $\pi$, and completes the proof
of the Proposition.
\qed

\begin{rem}
\label{rem_enddirectsum}
In the statement of Proposition \ref{prop_lowerweightcomp},
we can replace $\End(W)$ with $\End(W_1) \oplus \cdots \oplus \End(W_k)$,
and we can replace ${\id}_W$ with 
$c_1 {\id}_{W_1} \oplus \cdots \oplus c_k {\id}_{W_k}$, where 
$c_1, \dots, c_k$ are scalars. Indeed, when all $c_j = 1$,
this follows because $\End(W_1) \oplus \cdots \oplus \End(W_k)$ is
a $(G\times G)$-submodule of $\End(W)$, so we may compute the
closure in the smaller space. Moreover, in this embedding,
${\id}_W = {\id}_{W_1} \oplus \cdots  \oplus {\id}_{W_k}$. To see
that we can put scalars in front of each summand follows from
an easy analysis of the argument.
\end{rem}

\subsection{Lie algebra realization of the compactification}
\label{sec_liealgebra}

We give another realization of $X$, in which no choice of highest
weight is used. This realization is used in \cite{EL1, EL2} to
give a Poisson structure on the wonderful compactification.
 Let $n = \dim(G)$.
 $G\times G$ acts on $\Gr (n, \fg\oplus \fg)$ through
the adjoint action. Let $\fg_\Delta = \{ (x,x) : x\in \fg \}$,
the diagonal subalgebra. Then the stabilizer in $(G\times G)$
of $\fg_\Delta$ is $G_\Delta =\{ (g,g) : g\in G \}$, so
$(G\times G)\cdot \fg_\Delta \cong (G\times G)/G_\Delta \cong G$.

Let $\overline{G} = \overline{(G\times G)\cdot \fg_\Delta}$,
where the closure is computed in 
$\Gr (n,\fg\oplus \fg)$. Since $\Gr (n,\fg\oplus \fg)$ is
projective, $\overline{G}$ is projective.

\begin{prop}
\label{prop_liealgebrarealization}
$\overline{G} \cong X$. In particular, $\overline{G}$ is
smooth with $2^l$ orbits.
\end{prop}

To prove Proposition \ref{prop_liealgebrarealization}, we give
a representation theoretic interpretation of $\overline{G}$,
and then  apply the machinery developed above.

Embed $ i: \Gr (n,\fg\oplus \fg) \hookrightarrow 
\PR (\wedge^n (\fg \oplus \fg))$ via the Plucker embedding.
That is, if $U \in \Gr (n,\fg\oplus \fg)$ has basis
$u_1, \dots, u_n$, we map 
$U \to i(U)=[ u_1 \wedge \cdots \wedge u_n ]$. 
It is well-known that the Plucker embedding 
 is a closed embedding. 
Denote $[U]$ for $i(U) \in \PR (\wedge^n (\fg \oplus \fg))$.
Note that $\wedge^n(\fg \oplus \fg)$ is naturally a
representation of $G\times G$ by taking the exterior
power of the adjoint representation.

Clearly, $i:(G\times G)\cdot \fg_\Delta \to (G\times G)\cdot [\fg_\Delta]$
is an isomorphism. Since the Plucker embedding
is a closed embedding, it follows that
$i:\overline{(G\times G)\cdot \fg_\Delta} \to
 \overline{(G\times G)\cdot [\fg_\Delta]}$ is an isomorphism.

To prove Proposition \ref{prop_liealgebrarealization}, 
we apply Proposition \ref{prop_lowerweightcomp} to
prove $X\cong \overline{(G\times G)\cdot [\fg_\Delta]}$.

Choose  a nonzero vector  $v_\Delta$
in the line $[\fg_\Delta]$, and
let $E = U(\fg \oplus \fg)\cdot v_\Delta $.
 We show that as a $G \times G$-module, $E = \End(V(2\rho))
\bigoplus \oplus_i \End(V(\mu_i)) \bigoplus F$, where the
$V(\mu_i)$ are $G\times G$-modules such that $2\rho \ge \mu_i$,
 and $F$ is a $G\times G$-module.
Furthermore, we show that 
$[\fg_{\Delta}] = [{\id}_{V(2\rho)} + \oplus_i
 c_i {\id}_{V(\mu_i)} + 0]$, for some scalars $c_i$.
Then Proposition \ref{prop_lowerweightcomp} as refined
in Remark \ref{rem_enddirectsum} implies the Proposition.

To verify these assertions, we first analyze the $T\times T$-weights
in $\wedge^n(\fg\oplus \fg)$. Let $H_1, \dots, H_l$ be a basis
of $\ft$.
A basis of $T\times T$-weights on $\fg \oplus \fg$ is given by
the vectors:

\begin{enumerate}
\item $\{ (H_1,H_1), \dots, (H_l,H_l), (H_1,-H_1), \dots, (H_l,-H_l) \}$,
all of trivial weight $(0,0)$;

\item $\{ (E_\alpha , 0) : \alpha \in \Phi \}$, and note that $(E_\alpha, 0)$
generates the unique weight space weight $(\alpha,0)$;

\item $\{ (0, E_\alpha) : \alpha \in \Phi \}$, and note that $(0,E_\alpha)$
generates the unique weight space of weight $(0,\alpha)$.
\end{enumerate}

We obtain a basis of weight vectors of $\wedge^n(\fg\oplus \fg)$
as follows. Let $A$ and $B$ be subsets of $\Phi$ such that
$|A| + |B| \le n$ and let $R$ be a subset of the vectors in (1)
of cardinality $n - |A| - |B|$. For each triple $(A,B,R)$ as above,
we define a weight vector 
$$x_{A,B,R} = \wedge_{\alpha \in A} \ (E_\alpha,0) \bigwedge 
\wedge_{\beta \in B} \ (0,E_\beta)
\bigwedge \wedge_{i=1, \dots, n - |A| - |B| } \ K_i,$$
where the $K_i$ are the vectors in the subset $R$.
 Then the collection $\{ x_{A,B,R} \}$ is a basis
of weight vectors of $\wedge^n(\fg \oplus \fg)$, and it is routine
to check that the weight of $x_{A,B,R}$ is 
$(\sum_{\alpha \in A} \alpha , \sum_{\beta \in B} \beta)$.
 where
$A$ and $B$ are subsets of $\Phi$.

Take $A=\Phi^+$ and $B=-\Phi^+$ and 
let $R_0$ be the subset of basis vectors from (1) 
with  $K_1 = (H_1,H_1), \dots, K_l=(H_l,H_l)$. 
Then we let $v_0 = x_{\Phi^+, -\Phi^+, R_0}$, so $v_0$ is a weight
vector with weight $(2\rho, -2\rho)$. The weight of any $x_{A,B,R}$
is $(\lambda, \mu)$ where $\lambda \le 2\rho$ and $\mu \ge - 2\rho$.
It follows immediately that $v_0$ is a highest weight vector
of $\wedge^n(\fg\oplus \fg)$ of highest weight $(2\rho, -2\rho)$
relative to the Borel subgroup $B \times B^-$ of $G\times G$.
Further, note that
$[v_0]=[{\ft}_{\Delta} + \fu\oplus \fu^- ]$, where
${\ft}_{\Delta } =\{ (x,x): x\in \ft \}$, and
$\fu \oplus \fu^- $ is the Lie algebra of the unipotent radical
of $B \times B^-$.

Now we show that $v_0 \in E$. Choose $H \in \ft$ such that
$\alpha_i(H)=1$ for every simple root $\alpha_i$. Define
$\phi: \cstar \to T$ by $\phi(e^\zeta)=\exp (\zeta H)$, $\zeta
\in \C$. We claim that:

\begin{equation}
	\label{eq:triple_star}
	\lim_{z\to \infty} (\phi(z),e)\cdot [\fg_\Delta ] = [v_0] , 
\end{equation}
where $z\in \cstar$.

To check this claim, take 

 $$v_\Delta  = \wedge_{\alpha \in \Phi^+} (E_\alpha,E_\alpha) \bigwedge 
\wedge_{i=1, \dots, l } (H_i,H_i) \bigwedge
\wedge_{\beta \in \Phi^+} (E_{-\beta},E_{-\beta}),$$
and note that $[v_{\Delta}]=
[\fg_{\Delta}]$.

%Then $(\phi(z),e)v_\Delta$ is
%
%$$\wedge_{\alpha \in \Phi_+} (z^{k_\alpha} E_\alpha,E_\alpha) \bigwedge 
%\wedge_{i=1, \dots, l } (H_i,H_i) \bigwedge
%\wedge_{\beta \in \Phi^+} (E_{-\beta},E_{-\beta}),$$

%where $k_\alpha = \sum n_i$ if 
%$\alpha = \sum n_i \alpha_i$.

We decompose $(\phi(z),e) \cdot v_\Delta$
 into a sum of $2^{2r}$ terms ($2 r = \lvert \Phi \rvert$) 
 in the span of weight vectors $x_{A,B,R_0}$, with $R_0 
= \{ (H_1, H_1), \dots, (H_l, H_l) \}$.
To get the $2^{2r}$ terms, for each root $\gamma \in \Phi$,
decompose $(E_\gamma, E_\gamma)=(E_\gamma,0) + (0,E_\gamma)$,
and choose one of the two summands in each term in the product.
 Each term corresponds to a choice of the 
subset $A \subset \Phi$, and then $B = \Phi - A$.
For a root $\alpha$, let $ht(\alpha)= \sum k_i$, where
$\alpha = k_1 \alpha_1 + \dots + k_l \alpha_l$.
Then we compute
$$
(\phi(z),e)\cdot x_{A,B,R_0} = 
\wedge_{\alpha \in A} (z^{ht(\alpha)} E_\alpha, 0) \bigwedge 
\wedge_{i=1, \dots, l } (H_i,H_i) \bigwedge
\wedge_{\beta \in B} (0,E_{\beta}).$$

Thus $(\phi(z), e)\cdot x_{A,B,R_0} = z^{n_A} x_{A,B,R_0}$, where 
$n_A = \sum_{\alpha \in A} ht(\alpha)$. Let $n_0 = \sum_{\alpha \in \Phi^+} ht(\alpha)$. 
It follows easily from properties of roots that $n_0 \ge n_A$ for any $(A,B,R_0)$,
 and $n_0=n_A$ if and only
if $A=\Phi^+$ and $B=-\Phi^+$.

Then the formula 
$$
(\phi(z),e)\cdot [v_\Delta ] = [ z^{n_0} x_{\Phi^+ , -\Phi^+, R_0}
+ \sum_{A\not= \Phi^+ } z^{n_A} x_{A,B,R_0} ]$$
implies \eqref{eq:triple_star} by taking the limit as $z\to \infty$.

It follows easily that $v_0 \in E$. Indeed, $\PR (E)$ is
$T\times T$-stable and closed, so 
$\overline{(T\times T)\cdot [{\fg}_\Delta]} \subset \PR(E)$. 
It follows from \eqref{eq:triple_star} that $[v_0] \in \PR (E)$,
so $v_0 \in E$.

Now observe that since $v_0$ is a highest weight vector of
weight $(2\rho, -2\rho)$, 
$U(\fg\oplus \fg)\cdot v_0 \cong V(2\rho)\otimes V(-2\rho) \cong \End(V(2\rho))$
is a submodule of $E$. Moreover, by our determination of weights
of $\wedge^n(\fg\oplus \fg)$, any highest
weight vector occurring in $E$ has highest weight $\mu \le 2\rho$.
 It follows by separating out the irreducible $G\times G$
representations in $E$ of the form $\End(V(\mu))$ that

\begin{equation}
	\label{eq:quad_star}
	E = \End(V(2\rho)) \bigoplus \oplus_i \End(V(\mu_i)) \bigoplus F ,
\end{equation} 
where $F$ is a sum of irreducible
representation of $G\times G$ not isomorphic to $\End(V(\mu))$ for
any $\mu$.

It remains to verify the claim that
$[\fg_\Delta ] = [{\id}_{V(2\rho)} + \sum c_\mu {\id}_{V(\mu)} ]$.
Since $v_{\Delta}$ is $G_\Delta$-invariant, it follows that
its projection to each irreducible $(G\times G)$-representation
appearing in  \eqref{eq:quad_star} is $G_{\Delta}$-invariant. We recall the
classification of irreducible representations with a nonzero $G_\Delta$-fixed
vector.

\begin{lem}
\label{lem_diagonalinvariant}
Let $E$ be an irreducible representation of $G\times G$. Then $E$
has a nonzero $G_\Delta$-fixed vector if and only if $E\cong \End(W)$
for some irreducible representation $W$ of $G$. If $E = \End(W)$,
then its $G_\Delta$-fixed space  $\C \cdot {\id}_W$ is one dimensional.
\end{lem}

\noindent
{\bf Proof.}
Every irreducible representation of $G\times G$
is isomorphic to $V \otimes W^*$,
 where $V$ and $W$ are irreducible represenations
of $G$. Schur's Lemma implies that $V\otimes W^* \cong \Hom (W,V)$
has a nonzero $G_\Delta$-invariant vector if and only $W\cong V$,
so $V \otimes W^* \cong \End(V)$. Moreover, by Schur's Lemma, 
if $W = V$, the
space of invariants is generated by the identity.
\qed

In particular, by this Lemma, the projection of $v_\Delta$ to
each factor in \eqref{eq:quad_star} is a scalar times the identity.
Thus, we may write:

$[\fg_\Delta ] =  [c_\rho {\id}_{V(2\rho)} + \sum c_\mu {\id}_{V(\mu)} ]$.

To complete the proof, it suffices to show that $c_\rho$ is
nonzero. For this, it suffices to show 
the projection of $v_\Delta$ to $\End(V(2\rho))$
in \eqref{eq:quad_star} is nonzero. We verified in the proof of \eqref{eq:triple_star} that
when we write $v_\Delta$ as a sum of our standard basis vectors,
$v_0$ is a nonzero summand of weight $(2\rho, -2\rho)$.
The fact that the projection is nonzero 
follows using the $T\times T$-action and linear
independence of distinct eigenvalues. 
\qed

\begin{rem}
There is a unique $G\times G$-isomorphism
$\phi:X\to \overline{G}$ such that $[{\id}_V] \mapsto \fg_{\Delta}$.
We may use $\phi$ to find representatives for each of the $G\times G$-orbits
in $\overline{G}$.

The points $z_I$ in $X$ correspond to Lie subalgebras $\phi(z_I)=\fm_I$ 
of $\fg \oplus \fg$. Recall $\fl_I$, $\fu_I$ and $\fu_I^-$
from section \ref{sec_geometryoforbits}.
Then,
$$\fm_I = \{ (X + Z, Y + Z) : X \in \fu_I,
Y \in \fu_I^-, Z \in \fl_I \}.$$

This may verified as follows.
Find a curve $f:\PR^1 \to X $ with image in the toric variety
$X^\prime = \overline{(T\times \{ e \})\cdot [{\id}_V]}$
 such that $f(0)={\id}_V$ and $f(\infty)=z_I$. Using $T\times T$-equivariance,
consider the image of this curve in $\overline{G}$ and 
$\phi \circ f(\infty) = \phi(z_I)$. 
The morphism $\phi$ is completely determined by the $\phi(z_I)$ and 
$G\times G$-equivariance, and from this it is not difficult to
show that $\phi(z_I)=\fm_I$.

These subalgebras at infinity are used in
 the paper \cite{EL1}
to compute real points of the wonderful compactifiction for a particular
real form.
\end{rem}

\section{Cohomology of the compactification}
\label{sec_cohomology}

\subsection{$T\times T$-fixed points on $X$}
\label{sec_ttfixedpoints}

We explain how to compute the integral cohomology of $X$ following 
\cite{DCS2}.
We use the Bialynicki-Birula decomposition. Suppose that
 $Z$ is any smooth projective variety with a $\cstar$-action
and suppose the fixed point set $Z^{\cstar}$ is finite.
The idea is that we should be able to recover the topology
of $Z$ from the fixed point set $Z^{\cstar}$ and the action
of $\cstar$ on the tangent space at fixed points.

We make some remarks about the action. Suppose 
$m=\dim(Z)$. Let $z_0 \in Z$
be a fixed point for the $\cstar$-action. Then
$\cstar$ acts linearly on the tangent space $T_{z_0}(Z)$.
It follows that $T_{z_0}(Z)= \sum_{j=1}^m \C v_j $
where $a\cdot v_j = a^{n_j} v_j$. Then each $n_j$ is
nonzero, and the weights $(n_1, \dots, n_m)$ are
called the weights of $\cstar$ at $z_0$. 

Let $T_{z_0}^+(Z) = \sum_{n_j > 0} \C v_j$.

\begin{thm}
\label{thm_bbdecomposition} (Bialynicki-Birula \cite{BB}) \\
Let $Z$ be as above, and let $Z^{\cstar}=\{ z_1, \dots, z_n \}$.
For $i \in \{ 1, \dots, n\}$, let \\
$C_i = \{ z\in Z : \lim_{a\to 0} a\cdot z = z_i \}$. Then
\begin{enumerate}
	\item $Z = \cup_{i} \ C_i$

	\item $C_{z_i} \cong T_{z_i}^+(Z)$

	\item $C_{z_i}$ is locally closed in $Z$.
\end{enumerate}

\noindent As a consequence,
$H_*(Z) = \sum_{i=1}^n \Z \sigma_i,$ 
where $\sigma_i \in H_{2\dim (C_{z_i})}(Z)$
is the cycle corresponding to $\overline{C_{z_i}}$.
\end{thm}

For example, if $X=\PR (\C^{n+1})$, we may let $\cstar$ act on $X$
by 
$$a\cdot [(x_0, \dots, x_n)]=
[(a^n x_0, a^{n-1} x_1, \dots, a x_{n-1}, x_n )], a \in \cstar $$

Then the fixed points are the coordinate vectors and the corresponding
decomposition from Theorem \ref{thm_bbdecomposition} is the usual
cell decomposition of $\PR (\C^{n+1})$.

We want to do a similar analysis for $X = \overline{\psi(G)}$.
We will find a subgroup of $T\times T$ isomorphic to $\cstar$,
 prove that $X^{T\times T}=X^{\cstar}$, and compute the
Bialynicki-Birula decomposition. To do this, we need some
results on $X^{T\times T}$ and the action of $T\times T$
on the tangent space to the fixed points.

For $w\in W$, choose a representative $\dw \in N_G(T)$,
and for notational simplicity we use the identity of $G$
as a representative for the identity of $W$.
For simplicity, denote $z_0$ for $z_{\{ 1, \dots, l \}}$.
We remark that 
the unique closed orbit $(G\times G)\cdot z_0$
is isomorphic to  $G/B \times G/B^-$.

\begin{lem}
\label{lem_ttfixedpoints}
$X^{T\times T} = \cup_{(y,w) \in W\times W} \ z_{y,w}$, where
$z_{y,w}=(\dy,\dw)\cdot z_0$.
\end{lem}

\noindent
{\bf Proof.}
We first claim that $X^{T\times T} \subset (G\times G)\cdot z_0$.
Since $X = \cup_I \ (G\times G)\cdot z_I$, 
it suffices to consider any point
 $x\in (G\times G)\cdot z_I$. If $x\in X^{T\times T}$, 
then $(G\times G)_x$ must contain a torus of dimension $2l$.
But the stabilizer $(G\times G)_x$ is isomorphic to the
stabilizer $(G\times G)_{z_I}$. 
By the computation in Proposition \ref{prop_zIstabilizer},
a maximal torus of $(G\times G)_{z_I}$ is the set

$\{(t_1,t_2) \in T\times T : t_1 {t_2}^{-1} \in Z(L_I) \}$.

It is easy to see that this maximal torus has dimension
$l + |I|$. In particular, $(G\times G)_{z_I}$ contains
a $2l$-dimensional torus if and only if $I=\{ 1, \dots, l \}$.
The claim follows.

Recall the well-known fact that
$(G/B)^T = \{ \dw B : w\in W \}$. The Lemma now follows by
 using the isomorphism
$(G\times G)\cdot z_0 \cong G/B \times G/B^-$.
\qed

\begin{lem}
\label{lem_weightsatorigin}
Let $z_0 = z_{\{ 1, \dots, l\}}$. The weights of $T\times T$
on the tangent space $T_{z_0}(X)$ are 

(1) $(\frac{1}{\alpha}, 1) : \alpha \in \Phi^+ ;$

(2) $(1, \alpha) : \alpha \in \Phi^+ ;$

(3) $(\frac{1}{\alpha_i}, \alpha_i ): i= 1, \dots, l $

\end{lem}

\noindent
{\bf Proof.}
Since $X_0$ is an open neighborhood of $z_0$ in $X$,
the tangent space $T_{z_0}(X)=T_{z_0}(X_0)$. Since
$X_0 \cong U^- \times U \times Z$, 
$$T_{z_0}(X_0)\cong \fu^-\oplus \fu \oplus T_{z_0}(Z).$$
The first factor of $T$ acts on $\fu^-$ by the adjoint
action and the second factor acts trivially, and this
gives the weights appearing in (1). The first factor of $T$
acts trivially on $\fu$ and the second factor acts by
the adjoint action, which gives the weights appearing in (2).
We define a $T\times T$-action on $\C^l$ by
$$(t_1,t_2)\cdot (z_1,\dots ,z_l)=(\alpha_1(t_2/t_1)z_1, \dots ,
\alpha_l( t_2/t_1) z_l ), t_1, t_2 \in T.$$
We claim that $F:\C^l \to Z$ is $T\times T$-equivariant. 
The weights in (3) are easily derived from this claim.
The claim for $T\times \{ e \}$ follows the same
way as the derivation of the final
formula for $\psi(t)$ in the proof of Lemma \ref{lem_Zaffine}.
The claim for $\{ e \} \times T$ follows by computing
$(e,t)\cdot \psi(e)$ by the same method as in \ref{lem_Zaffine}. The
answer is different because $(e,t)\cdot v_k \otimes v_k^*
=v_k \otimes \frac{1}{\lambda_k(t)} v_k^*$, which is the
inverse of the corresponding formula for the $(t,e)$-action.
We leave the details to the reader.
\qed

Let $y, w \in W$.
Define $\phi :X \to X$ by $\phi (x)=(\dy, \dw)\cdot x$.
Then $z_{y,w}=\phi (z_0)$ is the $T\times T$-fixed
point of $X$ corresponding to $(yB,wB^-)$.

\begin{lem}
\label{lem_weightsatzyw}
 The weights of $T\times T$
on the tangent space $T_{z_{y,w}}(X)$ are 

(1) $(\frac{1}{y\alpha}, 1) : \alpha \in \Phi^+ ;$

(2) $(1, w\alpha) : \alpha \in \Phi^+ ;$

(3) $(\frac{1}{y\alpha_i}, w\alpha_i ): i= 1, \dots, l $

\end{lem} 

\noindent
{\bf Proof.}
Indeed, let $(t_1,t_2) \in T\times T$. It is routine to check that
$$
\phi((t_1,t_2)\cdot x)=(y(t_1), w(t_2))\phi(x).
$$
Thus, 
the differential $\phi_* :T_{z_0}(X)\to T_{z_{y,w}}(X)$ satisfies
the formula: 
$$\phi_* (t_1,t_2) = (y(t_1), w(t_2)) \phi_*.$$
In particular, if $\lambda_1, \dots, \lambda_n$ are the weights
of $T \times T$ at $z_0$, then $(y,w)\lambda_1, \dots, (y,w)\lambda_n$
are the $T\times T$-weights at $z_{y,w}$. The result now follows
from Lemma \ref{lem_weightsatorigin}.
\qed

\begin{rem}
\label{rem_toricclosure}
By Lemma \ref{lem_weightsatorigin},  it follows
 that the $\{ e \} \times T$
weights at the origin in $Z$ are $\alpha_i, i=1, \dots, l$.
By the proof of Lemma \ref{lem_weightsatzyw}, it follows that
the $\{ e \} \times T$ weights at the point $z_{w,w}$ are
$w(\alpha_i), i=1, \dots, l$.
Using these weight calculations, we may identify $\tilde{Z}:=
\cup_{w\in W} \ (\dw, \dw)\cdot Z$ with the closure of $Z$ in
$X$. Indeed, since $Z$ is smooth, it follows easily that $\tilde{Z}$
is a smooth toric variety for the torus $\{ e \} \times T \cong T$.
Every toric variety $Y$ with torus $S$ has a fan, which is a subset of 
${\fs}_{\R}=X_*(S)\otimes_{\Z} \R$. To compute the fan of a smooth
toric variety, we compute the weights $\lambda_1, \dots, \lambda_l$
of the action of the torus $S$ at each fixed point $z$. Define the
chamber $C_z$ for $z$ by $C_z = \{ x \in {\fs}_{\R} : \lambda_i(x)
\ge 0, i=1, \dots, l$\}. Then the fan of $Y$ is the subset
$F(Y):=\cup_{z\in Y^S} \ C_z$ together with the origin. The fan is called
complete if $F(Y)=\fs_{\R}$, and if $F(Y)$ is complete, then the
toric variety $Y$ is complete. For $\tilde{Z}$, the fan is easily
identified with the
Weyl chamber decomposition of ${\ft}_{\R}$, and in particular, it
is complete. It is easy to see that $Z$ is dense in $\tilde{Z}$,
and from this it follows that $\tilde{Z}$ coincides with the
closure of $Z$ in $X$. See \cite{fulton} for basic facts about
toric varieties.
\end{rem}

\subsection{Computation of cohomology}
\label{sec_computationofcohomology}

Recall the following well-known fact.

\begin{lem}
\label{lem_cstaractions}
Let $Y$ be a smooth $\cstar$-variety. Let $Y^\prime$ be a connected
component of $Y^{\cstar}$.
% and let $Y^\prime$
%be a connected component of $Y^{\cstar}$. 

(1) $Y^\prime$ is smooth and projective.

(2) $Y^\prime$ is a point if and only if there is a point $y \in
Y^\prime$ such that
 all the weights of $\cstar$ on $T_y(Y)$ are nontrivial.

\end{lem}

The condition that $\cstar$ has no trivial weights on $T_y(Y)$
is equivalent to the condition that the differentiated representation
of the Lie algebra $\C$ of $\cstar$ on $T_y(Y)$ has no zero weights.
It follows that to find a $\cstar$-action on $X$ with no fixed
points, it suffices to find some $A \in \ft\oplus \ft$ with
no zero weights on $T_{z_{y,w}}(X)$ for all pairs $(y,w) \in
W\times W$. The element $A$  generates $Lie(\cstar)=\C \cdot A$.
 Write $d\alpha \in \ft^*$ for the differential of
$\alpha$.

Choose $H\in \ft$ so that $d\alpha_i(H)=1$ for all simple roots
$\alpha_i$. It follows that if $\alpha$ is any root, $|d\alpha (H)| \ge 1$.
We may choose $n$ sufficiently large so that
$n > d\beta(H)$ for every root $\beta$. 
Then $| nd\alpha (H) | > |\beta(H)|$ for all roots $\alpha$ and
$\beta$. In particular,
\begin{equation}
\label{eq:ndH} 
	nd\alpha(H) > 0 \iff nd\alpha(H) + \beta(H) > 0
\end{equation}
for every pair of roots $\alpha$ and $\beta$.

 Consider $\eta: \cstar \to T\times T$ given by
 $\eta (a) =(\phi(a^n),\phi(a^{-1}))$. Then $d\eta (1)=(nH,-H)$.
 
We will let $\cstar$ act on the wonderful compactification
$X$ via $\eta$, so $a\cdot x = \eta(a) \cdot x$.

\begin{lem}
\label{lem_ttcstarfixedpoints}
$X^{T\times T} = X^{\cstar}$.
\end{lem}

\noindent
{\bf Proof.}
Indeed, it is clear that $X^{T\times T} \subset X^{\cstar}$.
Let $X^\prime$ be a connected component of $X^{\cstar}$.
Then $X^\prime$ is projective by Lemma \ref{lem_cstaractions} (1),
and is $T\times T$-stable since $T\times T$ commutes with
$\cstar$. The Borel fixed point theorem asserts that an
action of a connected solvable group on a projective variety
has a fixed point. Thus, $X^\prime$ contains a $T\times T$-fixed
point. By Lemma \ref{lem_ttfixedpoints}, some
$z_{y,w} \in X^\prime$.

We compute the weights of $(nH,-H)$ on the tangent space
$T_{z_{y,w}} (X)$ and show they are all nonzero. It follows
that all weights of $(nH,-H)$ on $T_{z_{y,w}}(X^\prime)$ are
nonzero. Then Lemma \ref{lem_cstaractions} (2) implies our claim.

By Lemma \ref{lem_weightsatzyw}, the eigenvalues of 
$(nH,-H)$ on $T_{z_{y,w}} (X)$ are:

(1) $-nd y\alpha (H), \alpha \in \Phi^+$

(2) $-d w\alpha (H), \alpha \in \Phi^+$

(3) $-(nd y\alpha_i(H) + d w\alpha_i(H)), \alpha_i$ simple.

By the choice of $H$, the numbers in (1) and (2) are nonzero.
By \eqref{eq:ndH} above, the numbers in (3) are nonzero.

\qed

For $y\in W$, let 
$$L(y)= | \{ \alpha_i \in \Delta : y\alpha_i \in -\Phi_+ \} |.$$

\begin{thm}
\label{thm_cohomologyofX}
$$H^*(X) = \sum_{(y,w)\in W\times W } \Z \sigma_{y,w},$$
where $\sigma_{y,w}$ has degree $2(l(y) + l(w) + L(y))$.
\end{thm}

\noindent
{\bf Proof.} 
For the  $\cstar$-action on $X$ from Lemma \ref{lem_ttcstarfixedpoints},
$X^{\cstar} = \{ z_{y,w} : (y,w) \in W \times W \}$. By Theorem
\ref{thm_bbdecomposition} and the remarks following Lemma 
\ref{lem_cstaractions}, it suffices to compute the number of
positive eigenvalues of the generator $(nH,-H)$ of the Lie
algebra of $\cstar$ on the tangent space at each fixed point.

By the definition of $H$, 
$-d y\alpha (H) > 0$ if and only if $y(\alpha) \in -\Phi^+$.
By the definition of $n$, $-(nd y\alpha_i(H) + d w\alpha_i(H)) > 0$
if and only if $y\alpha_i \in -\Phi^+$.
By the computation at the end of the proof of Lemma
\ref{lem_ttcstarfixedpoints}, it follows that the number of
positive eigenvalues of the generator $(nH,-H)$ at
the tangent space at $z_{y,w}$ is $l(y) + l(w) + L(y)$
(here we use the well-known fact that $l(y)$ is the number
of positive roots whose sign is changed by $y$). The Theorem
follows.
\qed

\begin{rem}
\label{rem_asymmetry}
The asymmetry of the roles of $y$ and $w$ in the statement of
the Theorem is due to the need to choose a $\cstar$-action
without fixed points.
\end{rem}

\section{Appendix on compactifications of general symmetric spaces}
\label{sec_appendix}

Let $G$ be complex semisimple with trivial center, and let
$\sigma:G \to G$ be an algebraic involution with fixed subgroup
$H=G^\sigma = \{ x\in G : \sigma(x)=x \}$. Then we call $(G,H)$
a symmetric pair and the associated homogeneous space $G/H$ is
called a semisimple symmetric space. As a special case, let $G_1$
be complex semisimple with trivial center and let $G=G_1 \times G_1$.
Then the diagonal subgroup $G_{1,\Delta}$ is the fixed point
 set of the involution $\sigma:G \to G$ given by
$\sigma(x,y)=(y,x)$. We refer to the symmetric pair
 $(G_1 \times G_1,G_{1,\Delta})$
as the group case, and note that the map 
$G_1 \cong (G_1 \times G_1)/G_{1,\Delta} \to G_1$, $(x,y) \mapsto 
xy^{-1}$ identifies the group $G_1$ as a symmetric space.
The purpose of this appendix is to explain a construction
of the DeConcini-Procesi compactification of the semisimple
symmetric space $G/H$. In the group case, this construction is
not identical to the construction given in Section 
\ref{sec_construction}, but it is equivalent.

First, we need some structure theory. For a $\sigma$-stable
maximal torus $T$, let 
$T^{-\sigma} = \{ t\in T : \sigma(t)=t^{-1}\}$.
 $T$ is called maximally split if $\dim(T^{-\sigma})$ is maximal
among all $\sigma$-stable maximal tori of $G$. All maximally split
tori are conjugate in $G$ \cite{springer}. For $T$ maximally split,
let $A$ be the connected component of the identity of
$T^{-\sigma}$ and let $r=\dim(A)$. Let $\fa$ be the Lie
algebra of $A$. We call $r$ the split rank
of the pair $(G,H)$. In the group case, if $T_1$ is a maximal
torus of $G_1$, $T_1 \times T_1$ is $\sigma$-stable and
$A =\{ (t,t^{-1}): t\in T_1 \}$, so the split rank is $\dim(T_1)$,
which is the rank of $G$.

Let $T$ be maximally split. Since $T$ is $\sigma$-stable, there
is an induced involution on the root system $\Phi$ of $(G,T)$, given
by $\sigma(\alpha)(t)=\alpha({\sigma}^{-1}(t))$.
Then $\sigma(\fg_\alpha)=\fg_{\sigma(\alpha)}$.
We say a root $\alpha$ is imaginary if $\sigma(\alpha)=\alpha$.
There exists a positive system $\Phi^+$ such that if
$\alpha \in \Phi^+$ and $\sigma(\alpha)\in \Phi^+$, then
$\alpha$ is imaginary (see [DCP], Lemma 1.2).
 Let $S = \{ \beta_1, \dots, \beta_s \}$ be the simple imaginary
roots for $\Phi^+$,
 and let $[S] = \{ \beta \in \Phi : \beta = \sum n_j \beta_j \}.$
Clearly, $\sigma$ acts as the identity on $[S]$.

\noindent CLAIM: If $\beta$ is an imaginary root, $\beta \in [S]$.

To prove the claim,
let $\alpha_1, \dots, \alpha_t$ be the simple nonimaginary roots
and define $f\in \ft^*=X(T)^*\otimes \C$ so that
 $f(\alpha_i)=1$ if
$\alpha_i$ is a simple non-imaginary root, and $f(\beta_j)=0$
if $\beta_j$ is imaginary. Note that if $\alpha$ is
a root, then 
\par\noindent (1)$f(\alpha) > 0 \iff \alpha \in \Phi^+ - [S]$,
\par\noindent (2)$f(\alpha)=0 \iff \alpha \in [S]$,  
\par\noindent (3)$f(\alpha) < 0 \iff
\alpha \in -\Phi^+ - [S]$.

Suppose $\beta$ is positive imaginary
and write 
$\beta = \sum_{n_i \ge 0} n_i \alpha_i + \sum_{m_j \ge 0} m_j \beta_j$, where the
first sum is over simple nonimaginary roots and the second sum is
over simple imaginary roots.
Then $\beta = \sigma(\beta)=\sum n_i \sigma(\alpha_i) + \sum m_j \beta_j$.
Since $\alpha_i$ is not imaginary, $\sigma(\alpha_i) \in -\Phi^+$.
Moreover, $\alpha_i \not\in [S]$, so $\sigma(\alpha_i)\not\in [S]$.
It follows that $f(\sigma(\alpha_i)) < 0$, and hence that $f(\sigma(\beta))
< 0$ if some $n_i \not= 0$. But $f(\beta) \ge 0$, so $f(\beta)=0$
and
all $n_i = 0$, so $\beta \in [S]$.

Let $\fm = \ft + \sum_{\alpha \in [S]} \fg_{\alpha}$. Then $\fm$ is
a Levi subalgebra whose roots form the full set of imaginary
roots.
Moreover, $\fm $ is the centralizer
of $\fa$.
Let $\fp = \fm + \fu$, where $\fu$ is the unipotent radical
of the Borel subalgebra determined by $\Phi^+$. Let $\fn$
be the nilradical of $\fp$, and let $\fn^-$ be the opposite $T$-stable
nilradical, so $\fg = \fp + \fn^-$. Then 
$\fn$ (resp. $\fn^-$) is the sum of the positive (resp. negative) 
 nonimaginary
root spaces.  Let $P, M, N$ and $N^-$ be the corresponding connected
groups. If $\fg_{\R}$ is the real form
of $\fg$ such that $\fh$ is the complexification of a
maximal compact subalgebra of $\fg_{\R}$, then $\fp$ is
the complexification of a minimal parabolic subalgebra of
$\fg_{\R}$, and $\fa$ and $\fn$ are complexifications of
corresponding factors in an Iwasawa decomposition of
$\fg_{\R}$.

 Let $\alpha_i$ be a non-imaginary
simple root. Then $\sigma(\alpha_i)=-\alpha_j - \sum_{i=1}^s k_i \beta_i$
for some nonimaginary simple root $\alpha_j$ and some 
nonnegative integers $k_i$ (see [DCP], Lemma 1.4).
 In this case, we define 
$\sigma(i)=j$, and note that $\sigma(j)=i$.
 $\sigma$ defines an involution on the indexing set of nonimaginary
simple roots 
$\{ 1, \dots, t \}$ with $r$ orbits. Renumber the 
$\alpha_1, \dots, \alpha_t$ so that each orbit of $\sigma$ on
$\{ 1, \dots, t \}$
has one element in $\{ 1,\dots, r\}$.
 If $\beta$ is imaginary, then
$d\beta$ vanishes on $\fa$, and
   $\sigma(d\alpha_i)(Y) = 
-d\alpha_j (Y)$ for all $Y\in \fa$. By definition, the restricted roots
 are the set of nontrivial characters 
$\{ \alpha \lvert_{A} : \alpha \in \Phi \}$. 
For $\alpha$ non-imaginary, set $\overline{\alpha}=\alpha\lvert_{A}$.
We call $\overline{\alpha_1}, \dots, \overline{\alpha_r}$ the simple
restricted roots.
Then each restricted root is a sum of simple restricted roots.
We may use the direct sum decomposition $\ft = \ft^\sigma + \fa$
to identify $\fa^* \cong \{ f \in \ft^* : f(\ft^\sigma)=0 \}$.
For $\alpha$ nonimaginary, $d\overline{\alpha}$ corresponds
to $\frac{d\alpha - \sigma(d\alpha)}{2}$ in this identification.
If $\sigma(\alpha_i)=-\alpha_i$,
 we call $\alpha_i$ a real root, and if $\sigma(\alpha_i)
\not= \pm \alpha_i$, we call $\alpha_i$ a complex root. The Satake diagram
of $(\fg,\sigma)$ is the Dynkin diagram of $\fg$ with simple
imaginary roots colored black, and other simple roots colored
white, with a double-edged arrow connecting $\alpha_i$ and $\alpha_j$ if
$\sigma(i)=j$. The Satake diagram determines $\sigma$
up to isomorphism \cite{araki}.

Let $\tilde{H}$ be the preimage of $H$ in the simply connected
cover $\tilde{G}$ of $G$. Let $V$ be an irreducible representation
of $\tilde{G}$ with highest weight vector $v_0$ of weight
$\lambda$. The induced $\tilde{G}$-action on $\PR (V)$ 
factors to give a $G$-action on $\PR (V)$ and similarly, we obtain
a $G\times G$-action on $\PR (\End(V))$.
 We choose $\lambda$ so that $G_{[v_0]}=P$,
where $G_{[v_0]}$ is the stabilizer in $G$ of the  line
$[v_0]$ through $v_0$. Equivalently, $\lambda(\alpha_i)>0$
for every simple nonimaginary root $\alpha_i$, and $\lambda(\beta_i)
= 0$ for every simple imaginary root $\beta_i$. In \cite{DCP},
such a weight $\lambda$ is called {\it regular special}.
Complete $v_0$ to a basis $v_0, \dots, v_n$ of weight vectors,
and number the weights so $v_i$ has weight $\lambda - \alpha_i$,
for $i=1, \dots, r$. Choose also a dual basis $v_0^*, v_1^*,
\dots, v_n^*$. Let $\PVopen = \{ [v] \in \PR (V) : v_0^* (v)
\not= 0$\}. The hypothesis on $\lambda$ implies that the morphism
$N^- \to N^- \cdot [v_0]$ is an isomorphism, and the image is
closed in $\PVopen$ by Lemma \ref{lem_unipotentorbitclosed}. 

Consider the twisted conjugation action of $G$ on
$\PR (\End(V))$ given by $g\cdot [C] = g [C] \sigma(g^{-1})$
for $g\in G$ and $C\in \End(V)$. Let
$h={\id}_V$. Then the stabilizer $G_{[h]}$ of the
line $[h]$ is $H$, so the orbit $G\cdot [h] \cong G/H$,
and we may regard the irreducible subvariety $X:= \overline{G\cdot [h]}$ as
a compactification of $G/H$. We wish to show $X$ is smooth
and describe its geometry.

For this, let $\openp = \{ [C] \in \PR : v_0^*(C \cdot v_0) \not= 0 \}$
and let $X_0 = X \cap \openp$. As in
\ref{sec_structureopenaffine}, we identify $X_0$ with an affine space,
and then show $X$ is a union of translates of $X_0$. For this,
note that the identification $V\otimes V^* \cong \End(V)$ is
$G$-equivariant, where we use the twisted conjugation action
on $\End(V)$ and the $G$-action on $V\otimes V^*$ is the
unique action such that
 $g\cdot (v\otimes f)=g\cdot v \otimes \sigma(g)\cdot f$
(the $G$-action on $V^*$ is the usual action on the dual).
We now compute the $A$-action on $[h]$. We identify
${\id}_V = \sum v_i \otimes v_i^*$ and obtain
$$
t\cdot [h] = t\cdot [\sum v_i \otimes v_i^*] = [v_0 \otimes v_0^* +
\sum_{i=1}^n \frac{1}{(\gamma_i - \sigma(\gamma_i))(t)} v_i \otimes v_i^* ],$$
where $\lambda - \gamma_i$ is the weight of $v_i$,
by reasoning as in \ref{lem_Zaffine}. For $j > 0$, 
$\gamma_j - \sigma(\gamma_j) = 
\sum_{i = 1}^r n_{ij } (\alpha_i - \sigma(\alpha_i))$ for
nonnegative integers $n_{ij}$.
It follows that the morphism $\psi:\C^r \to \openp$
given by
$$
(z_1 , \dots, z_r) \mapsto [v_0 \otimes v_0^* +
 \sum_{i=1}^r z_i v_i \otimes v_i^* + 
\sum_{j=r+1}^n \prod_{i=1}^r {z_i}^{n_{ij}} v_j \otimes v_j^* ]
$$
is an isomorphism. Let $Z=\psi(\C^r)$. Note that $A$ acts on $Z$
through its quotient $A/D$, where $D=\{ a\in A : \alpha_i(a)=\pm 1,
i=1, \dots, r\}$.

Note that $N^-$ preserves $\openp$ so $N^-$ acts on $X_0$. Indeed,
if $v_0^*(C\cdot v_0)\not= 0$, then for $n\in N^-$,
 $v_0^*(n\circ C \circ \sigma(n)^{-1}\cdot v_0)=v_0^*(n \cdot C\cdot v_0)
=v_0^* ( C\cdot v_0)$
since $\sigma(n)\in N$ and $n$ fixes the lowest weight vector
$v_0^*$ of $V^*$. Thus, we have the morphism $\chi:N^- \times \C^r
\to X_0$ given by $\chi(u,z)=u\cdot \psi(z)$. We claim that
$\chi$ is an isomorphism. 
For this, by reasoning as in the proof of 
\ref{thm_X0identifier},
 it suffices to construct a $N^-$-equivariant morphism
$\beta:X_0 \to N^-$ such that $\beta \circ \chi (u,z)=u$.
To construct $\beta$, consider the map $\nu:\openp \to \PVopen$
defined by $\nu([C])=[C\cdot v_0]$. Note that
$$
\nu(N^-A \cdot [h] )=\{ \nu([x\sigma(x)^{-1}]) : x\in N^- A \} =
$$
$$
\{ [ x\sigma(x)^{-1}\cdot v_0 ] : x\in N^- A \} = 
\{ [x\cdot v_0 ] : x\in N^- A \} = N^- \cdot [v_0].
$$
It is easy to check that $\dim(N^-A \cdot [h]) = \dim G\cdot [h]$,
and it follows that $N^-A \cdot [h] $ is open and dense in $X_0$,
since an orbit of an algebraic group on a variety is open in its
closure.
It follows that $\nu(X_0)\subset N^- \cdot [v_0]$, since
$N^- \cdot [v_0]$ is closed in $\PVopen$. We thus obtain 
a morphism $\beta: X_0 \to N^-$ by composing $\nu$ with
the isomorphism $N^- \cdot [v_0] \to N^-$. It is easy to 
check that $\beta$ has the required property, and this gives
the following theorem.

\begin{thm}
\label{thm_X0symmetric}
$\chi:N^- \times \C^r \to X_0$ is an isomorphism of varieties.
In particular, $X_0$ is smooth.
\end{thm}

To prove $X$ is smooth, let $V = \cup_{g\in G} \ g\cdot X_0$.
We want to show $V=X$. For this, the following result from
\cite{DCS1} is useful.

\begin{prop}
\label{prop_showcomplete}

Let $X$ be a complex projective
$G$-variety. Let $x\in X$ be a point with stabilizer $G_x = H$.
Let $Y \subset X$ be a locally closed $G$-stable subset containing the orbit
$G\cdot x$ and suppose that

(a) $G\cdot x$ is dense in $Y$

(b) The closure of $T\cdot x$ in $Y$ is projective.

Then $Y$ is projective.
\end{prop}

The proof, given in \cite{DCS1}, uses the valuative criterion for
properness and a variant of the Hilbert-Mumford criterion.

To apply this result to show $V=X$, it suffices to show that the
closure $X^\prime$ of $T\cdot [h] $ in $V$ is complete. For this, we use
the theory of toric varieties. Let $W_A = N_G(A)/Z_G(A)$, the 
so-called little Weyl group. For $w\in W_A$, let $\dot{w}$ be
a representative in $N_G(A)$. Clearly $Z \subset X^\prime$, and
it follows that $\cup_{w\in W_A} \ \dot{w}\cdot Z$ is contained
in $X^\prime$.  But $\cup_{w\in W_A} \ \dot{w}\cdot Z$ is a smooth toric
variety for the torus $A/D$, and it is straightforward 
as in Remark \ref{rem_toricclosure}
   to show that its fan is the
Weyl chamber decomposition of $\fa_{\R}$. In particular, the fan
is complete, so  $\cup_{w\in W_A} \ \dot{w}\cdot Z$ is a complete
variety by standard results on toric varieties \cite{fulton}.
Since $\cup_{w\in W_A} \ \dot{w}\cdot Z$ is dense in $X^\prime$,
it folows easily that $X^\prime$ is complete, so it is projective.
Now Proposition \ref{prop_showcomplete} implies that $V$ is projective,
so $V=X$.

Further, let $Z_i = \psi(\{(z_1, \dots, z_r): z_i \not= 0 \}$,
and let $D_i = \overline{G\cdot Z_i}$. Then by arguments 
as in \ref{sec_orbitdescription}, each $D_i$ is a smooth divisor,
and the different $D_i$ meet transversally.
 For $I \subset \{ 1, \dots, r\}$, let $D_I = \cap_{i\in I} D_i$.

\begin{thm}
\label{thm_globalgeometrysymmetric} The following hold:
\begin{enumerate}
	\item $X$ is smooth.

	\item Every $G$-orbit closure in $X$ is $D_I$ for
	some $I$. In particular, all $G$-orbit closures are smooth.

	\item For each $G$-orbit closure $D_I$ in $X$, there is a parabolic 
	subgroup $Q_I$ such that $D_I$ fibers over $G/Q_I$ with
	fiber a wonderful compactification of a symmetric space
	of the adjoint quotient of the  Levi factor of $Q_I$.
\end{enumerate}
\end{thm}

%%%%%
%
% BIBLIOGRAPHY
%
%%%%%

\bibliographystyle{bfjthesis}
\bibliography{dpbib}

\end{document}